\documentstyle[11pt]{article}
\input{amssym.def}
\input{amssym}
\newcommand\C{{\mbox{$\Bbb C$}}}
\newcommand\R{{\mbox{$\Bbb R$}}}

\newcommand\qed{\begin{flushright} $\Box$ \end{flushright}}

\newtheorem{pr}{Proposition}
\newtheorem{de}{Definition}
\newtheorem{th}{Theorem}

\newcommand{\be}{\begin{equation}}
\newcommand{\ee}{\end{equation}}
\newcommand{\bea}{\begin{eqnarray}}
\newcommand{\eea}{\end{eqnarray}}
\newcommand{\bean}{\begin{eqnarray*}}
\newcommand{\eean}{\end{eqnarray*}}
\newcommand{\var}{\varepsilon}
\newcommand{\suml}{\sum\limits}
\newcommand{\ecke}{\;_-\!\rule{0.2mm}{0.2cm}\;\;}
\begin{document}
\begin{LARGE}
\begin{center}
{\bf Lorentzian twistor spinors and CR-geometry}
\end{center}
\end{LARGE}
\vspace{0.5cm}
Helga Baum\\
\begin{small} {\it Humboldt University, Department of Mathematics, Sitz:
Ziegelstrasse 13a, 10099 Berlin, Germany}\\[0.2cm]
{\it Abstract}:
We prove that there exist global solutions of the twistor equation on
the Fefferman spaces of strictly pseudoconvex spin manifolds of
arbitrary dimension and we study their properties.\\[0.2cm]
{\it Keywords}: Twistor equation, Twistor spinors, Lorentzian manifolds, CR-geometry,
Fefferman spaces.\\[0.2cm]
{\it MS classification}: 58G30, 53C50, 53A50.
\end{small}

\section{Introduction} \label{s1}

In the present paper we study a relation between the Lorentzian
twistor equation and CR-geometry.
Besides the Dirac operator there is a second important conformally
covariant differential operator acting on the spinor fields
$\Gamma(S)$ of a smooth
semi-Riemannian spin manifold $(M,g)$ of dimension $n$ and index $k$,
the so-called {\em twistor operator} ${\cal D}$. The twistor operator
is defined as the composition of the spinor
derivative $\nabla^S$ with the projection $p$ onto the kernel of the
Clifford multiplication $\mu$
\[ {\cal D}:\Gamma(S)\stackrel{\nabla^S}{\longrightarrow}
\Gamma(T^*M\otimes S) \stackrel{g}{\approx} \Gamma(TM\otimes
S)\stackrel{p}{\longrightarrow} \Gamma(\mbox{Ker}\,\mu). \]
The elements of the kernel of ${\cal D}$ are called {\em twistor
spinors}.
A spinor field $\varphi$
is a twistor spinor if and only if it satisfies the
{\em twistor equation} \[ \nabla^S_X \varphi + \frac{1}{n} X \cdot
D \varphi = 0 \] for each vector field $X$, where
$D$ is the Dirac operator. Each twistor spinor $\varphi$ defines a
conformal vector field $V_{\varphi}$ on $M$ by
\[ g(V_\varphi , X) = i^{k+1}\,\langle X \cdot \varphi , \varphi
\rangle \,. \]
Twistor spinors were introduced by R.Penrose in General Relativity
(see \cite{Penrose:67}, \cite{Penrose/Rindler:86},
\cite{Nieuwenhuizen/Warner:84}). They are related to Killing vector
fields in semi-Riemannian supergeometry (see
\cite{Alekseevski/Cortes/ua:97}).
In Riemannian geometry the twistor equation first appeared as an
integrability condition for the canonical almost complex structure of
the twistor space of an oriented four-dimensional Riemannian manifold
(see \cite{Atiyah/Hitchin/ua:78}).
In the second half of the 80th Lichnerowicz and
Friedrich started the systematic investigation of twistor spinors on
Riemannian spin manifolds from the view point of conformal
differential geometry. Nowadays one has a lot of structure results and
examples for ma\-ni\-folds with twistor spinors in the Riemannian setting
(see \cite{Lichnerowicz1:88}, \cite{Lichnerowicz2:88},
\cite{Lichnerowicz1:89}, \cite{Friedrich:89} \cite{Lichnerowicz2:90},
\cite{Friedrich/Pokorna:91}, \cite{Baum/Friedrich/ua:91},
\cite{Habermann:90}, \cite{Habermann:93}, \cite{Habermann:94},
\cite{Kuehnel/Rademacher:94},
\cite{Kuehnel/Rademacher:95}, \cite{Kuehnel/Rademacher1:96},
\cite{Kuehnel/Rademacher2:96}). Crucial results were obtained by
studying the properties of the conformal vector field
$V_{\varphi}$ of a twistor spinor $\varphi$. Twistor operators also
turned out to be a usefull tool in proving sharp eigenvalue estimates
for coupled Dirac operators on compact Riemannian manifolds (see eg
\cite{Baum1:94}). \\
In opposite to this, there is not much known
about solutions of the twistor equation in the general Lorentzian
setting. In 1991 Lewandowski studied local solutions of the twistor
equation on 4-dimensional space-times, (\cite{Lewandowski:91}). In
particular, he proved that a 4-dimensional space-time admitting a
twistor spinor $\varphi$ without zeros and with twisting conformal
vector field $V_{\varphi}$ is locally conformal equivalent to a
Fefferman space. On the other hand, on 4-dimensional Fefferman spaces
there exist {\em local} solutions of the twistor equation.
The aim of the present paper is the generalisation of this result. \\
Fefferman spaces were defined by Fefferman (\cite{Fefferman:76}) in case
of strictly pseudoconvex hypersurfaces in $\C^n$, its definition was
extended by Burns, Diederich, Shnider (\cite{Burns/Diederich/ua:77}),
Farris (\cite{Farris:86}) and Lee (\cite{Lee:86}) to general
non-degenerate CR-manifolds. Sparling (\cite{Sparling:85}), Lee
(\cite{Lee:86}), Graham (\cite{Graham:87}) and Koch (\cite{Koch:88})
studied geometric properties of Fefferman spaces. A Fefferman space is
the total space of a certain $S^1$-principal bundle over a
non-degenerate CR-manifold $M$ equipped with a semi-Riemannian metric
defined by means of the Webster connection.
By changing the topological type of the $S^1$-bundle defining the
Fefferman space, we can prove that there are {\em global} solutions of the
twistor equation on the (modified) Fefferman spaces of strictly
pseudoconvex spin manifolds of arbitrary dimension. These solutions
have very special geometric properties which are only possible on
Fefferman spaces. More exactly, we prove (see Theorem \ref{t1},
Theorem \ref{t2}):\\[0.3cm]
{\em Let $\,(M^{2n+1},T_{10},\theta)\,$ be a strictly pseudoconvex
spin manifold and $\,(\sqrt{F},h_\theta)\,$ its Fefferman space. Then,
on the Lorentzian spin manifold $\,(\sqrt{F},h_\theta)\,$ there exist
a non-trivial twistor spinor $\phi$ such that
\begin{enumerate}
\item  The canonical vector field $V_\phi$ of $\phi$ is a regular
isotropic Killing vector field.
\item $\, V_\phi \cdot \phi = 0\,$. In particular, $\phi$ is a pure or
partially pure spinor field.
\item $\,\nabla_{V_\phi} \phi = i\,c\,\phi\,, \quad c=
\,\mbox{const}\,\in \R
\setminus \{0\}\,$.
\end{enumerate}
On the other hand, if $(B,h)$ is a Lorentzian spin manifold with a
non-trivial twistor spinor satisfying 1. - 3.,
then $B$ is an $S^1$-principal bundle over a stricly
pseudoconvex spin manifold $\,(M,T_{10},\theta)\,$ and $(B,h)$ is
locally isometric to the Fefferman space $\,(\sqrt{F},h_\theta)\,$ of
$(M,T_{10},\theta)\,$}.\\[0.3cm]
In particular, if $(M^{2n+1},T_{10},\theta)$ is a compact strictly
pseudoconvex spin manifold of constant Webster scalar curvature, then
the Fefferman space $(\sqrt{F},h_\theta)$ of $(M,T_{10},\theta)$ is a
(2n+2)-dimensional non-Einsteinian
Lorentzian spin manifold of constant scalar curvature $R$ and the
twistor spinor $\phi$ defines eigenspinors of the Dirac operator of
$(\sqrt{F},h_\theta)$ to the eigenvalues $\, \,\pm
\frac{1}{2}\sqrt{\frac{2n+2}{2n+1}R}\,\,$  with constant length.\\

After some algebraic prelimeries in section 2 we introduce in section
3 the notion of Lorentzian twistor spinors and explain some of their
basic properties. In order to define the (modified) Fefferman space we
recall in section 4 the basic notions of pseudo-hermitian geometry. In
particular, we explain the properties of the Webster connection of a
non-degenerate pseudo-hermitian manifold, which are important for
the spinor calculus on Fefferman spaces. In section 5 the Fefferman
spaces are defined and in section 6 we derive a spinor calculus for
Lorentzian metrics on $S^1$-principal bundles with isotropic fibre
over strictly pseudoconvex spin manifolds.
Finally, section 7 contains the proof of the Theorems 1 and 2 which
state the properties of the solutions of the twistor equation on
Fefferman spaces of strictly pseudoconvex spin manifolds.

\section{Algebraic prelimeries} \label{s2}

For  concrete calculations we will use the following
realization of the spinor representation.
Let $\,\mbox{Cliff}_{n,k}\,$ be the Clifford algebra of
$\,(\R^n,-\langle \cdot, \cdot \rangle_k)\,$, where $\,\langle \cdot, \cdot
\rangle_k\,$ is the scalar product
$\,\, \langle
x,y\rangle_k:=-x_1y_1-\ldots-x_ky_k+x_{k+1}y_{k+1}+\ldots+x_ny_n\,\,$.
For the canonical basis
$(e_1,\ldots,e_n)$ of $\R^n$ one has the following relations in
$\,\mbox{Cliff}_{n,k}\,:\,\,
 e_i\cdot e_j+e_j\cdot
e_i=-2\var_j\delta_{ij},\,\,$  where $ \var_j=\left\{\begin{array}{rl}
-1 & j\le k\\ 1 & j> k \,\end{array}\right. $.
Denote $\,\,\tau_j = \left\{\begin{array}{ll} i & j\le k\\ 1 &
j>k\end{array}\right.\,\,$ and
\begin{displaymath}
U = \left(\begin{array}{cc} i & 0\\ 0 & -i\end{array}\right),\quad
V=\left(\begin{array}{cc} 0 & i\\ i & 0\end{array}\right),\quad
E=\left(\begin{array}{cc} 1 & 0\\ 0 & 1\end{array}\right),\quad
T=\left(\begin{array}{cc} 0 & -i\\ i &
0\end{array}\right).
\end{displaymath}
Then an isomorphism \[ \phi_{2m,k}:\mbox{Cliff}^{\Bbb C}_{2m,k}
\longrightarrow M(2^m;\C) \] is given by the Kronecker product
\begin{eqnarray}\label{1}
\begin{array}{llll} \phi_{2m,k}(e_{2j-1}) & = &
\tau_{2j-1} &E\otimes\ldots\otimes E \otimes U \otimes
T\otimes\ldots\otimes T\\ \phi_{2m,k}(e_{2j}) & =
& \tau_{2j} &E\otimes\ldots\otimes E\otimes V \otimes
\underbrace{T\otimes\ldots\otimes T}_{j-1} \end{array}.
\end{eqnarray}
Let
$\mbox{Spin}_0(n,k)\subset \mbox{Cliff}_{n,k}$ be the connected
component of the identity of the spin group. The spinor representation is given
by \[ x_{n,k}=\phi_{n,k}|_{\mbox{Spin}_0(n,k)}:\mbox{Spin}_0(n,k)
\longrightarrow \mbox{GL} (\C^{2^m}). \]
We denote this representation by
$\,\Delta_{n,k}\,$. If $n=2m$, $\,\Delta_{2m,k}\,$ splits into the sum
$\,\Delta_{2m,k}=\Delta^+_{2m,k}\oplus \Delta^-_{2m,k}\,$, where
$\,\Delta^\pm_{2m,k}\,$ are the eigenspaces of the endomorphism
$\,\phi_{2m,k}(e_1\cdot\ldots\cdot e_{2m})\,$ to the eigenvalue $\pm
i^{m+k}$.
Let us denote by $u(\delta)\in\C^2$ the vector $\,
u(\delta)=\frac{1}{\sqrt{2}}{1\choose -\delta i},\,\,\delta=\pm 1
\,\,$ and let \be\label{2}
u(\delta_1,\ldots,\delta_m)=u(\delta_1)\otimes\ldots\otimes
u(\delta_m)\qquad \delta_j=\pm 1. \ee Then
$\;(u(\delta_1,\ldots,\delta_m)\,|\,\prod\limits^m_{j=1}\delta_j=\pm
1)\;$ is an orthonormal basis of $\,\Delta^\pm_{2m,k}\,$ with respect
to the standard scalar product of $\C^{2^m}$.

\section{Lorentzian twistor spinors} \label{s3}

Let $(M^{n,1},g)$ be a connected space- and time oriented Lorentzian
spin manifold with a fixed time orientation $\,\xi\in \Gamma(TM)\,$,
$g(\xi,\xi)=-1$. We denote by $S$ the spinor bundle of $(M^{n,1},g)$,
by $\,\nabla^S:\Gamma(S) \to \Gamma(TM^*\otimes S)\,$ the spinor
derivative given by the Levi-Civita connection of $(M^{n,1},g)$ and by
$\,D:\Gamma(S)\to\Gamma(S)\,$ the Dirac operator on $S$.\\ On $S$ there
exists an indefinite scalar product $\,\langle\cdot,\cdot\rangle\,$ of
index $\,\frac{1}{2}\dim S\,$ such that \bea\label{3} \langle
X\cdot\varphi,\psi\rangle &=& \langle\varphi, X\cdot\psi\rangle\\
X\langle\varphi,\psi\rangle &=&
\langle\nabla^S_X\varphi,\psi\rangle+
\langle\varphi,\nabla^S_X\psi\rangle \label{4} \eea for all vector
fields $X$ and all spinor fields $\varphi,\psi
\in\Gamma(S)$. Furthermore, there is a positive definite scalar
product $\,(\cdot ,\cdot )_\xi\,$ on $S$ depending on the time
orientation $\xi$ such that \be\label{5} \langle\varphi,\psi\rangle
=(\xi\cdot\varphi,\psi)_\xi \ee for all $\varphi,\psi\in\Gamma(S)$
(see \cite{Baum:81}, chap.1.5, 3.3.1.). Let $\,p:TM\otimes S
\longrightarrow \mbox{Ker}\,\mu\,$ denote the orthogonal projection onto
the kernel of the Clifford multiplication $\mu$ (with respect to
$\langle \cdot , \cdot\rangle$). $p$ is given by \[
p(X \otimes \varphi)=X \otimes \varphi+\frac{1}{n}\suml^n_{k=1}\var_k
s_k \otimes s_k \cdot X \cdot \varphi, \] where $(s_1,\ldots,s_n)$ is
a orthonormal
basis
of $(M,g)$ and $\var_k=g(s_k, s_k)=\pm 1$.

\begin{de}
The twistor operator ${\cal D}$ of
$\,(M^{n,1}, g)\,$ is the operator given by the composition of the
spinor
derivative with the projection $p$ \[ {\cal
D}:\Gamma(S)\stackrel{\nabla^S}{\longrightarrow} \Gamma(T^*M\otimes
S) \stackrel{g}{\approx} \Gamma(TM\otimes S)\stackrel{p}{\longrightarrow}
\Gamma(\mbox{Ker}\,\mu). \]
\end{de}
Locally, we have \[ {\cal
D}\varphi=\suml^n_{k=1}\var_k s_k\otimes(\nabla^S_{s_k}\varphi+
\frac{1}{n}s_k\cdot D\varphi). \]

\begin{de}
A spinor field $\varphi\in\Gamma(S)$ is called a
twistor spinor, if ${\cal D}\varphi=0$.
\end{de}
Let us first recall some properties of twistor spinors which are
proved in the same way as in the Riemannian case.

\begin{pr} \label{pr1} {\em (\cite{Baum/Friedrich/ua:91}, Th.1.2)}\\
For a spinor field  $\varphi\in
\Gamma(S)$ the following conditions are equivalent:\\[0.1cm]
\begin{tabular}{cl}
1. & $\varphi$ is a twistor spinor.\\[0.1cm]
2. & $\varphi$ satisfies
the so-called twistor equation
\end{tabular}
\vspace{-0.05cm}\\
\begin{equation} \label{6}
\nabla^S_X\varphi+\frac{1}{n}X\cdot D\varphi=0
\end{equation}
\vspace{-0.4cm}\\
\begin{tabular}{cl}
   & for all vector fields $X$.\\[0.1cm]
3. & For all
vector fields $X$ and $Y$
\end{tabular}
\vspace{-0.05cm}\\
\be\label{7}
X\cdot\nabla^S_Y\varphi+Y\cdot\nabla^S_X\varphi=\frac{2}{n}\,g(X,Y)\,
D\varphi \ee
\vspace{-0.4cm}\\
\begin{tabular}{cl}
   & holds. \\[0.1cm]
4. & There exists a spinor field $\psi\in\Gamma(S)$ such that
\end{tabular}
\vspace{-0.05cm}\\
\be\label{8}
\psi=g(X,X)X\cdot\nabla^S_X\varphi \ee
\vspace{-0.4cm}\\
\begin{tabular}{cl}
   & $\;\;\;$for all vector fields $X$ with $\,|g(X,X)|=1$.
\end{tabular}
\end{pr}
\vspace{0.1cm}

\begin{pr} \label{pr2} {\em (\cite{Baum/Friedrich/ua:91}, Th.1.7)}\\
The twistor operator is
conformally covariant: Let $\,\tilde{g} = e^{2\sigma} g\,$ be a
conformally equivalent metric to $g$ and let $\tilde{D}$ be the
twistor operator of $(M,\tilde{g})$. Then \[ \tilde{D} \tilde{\varphi}
= e^{- \frac{1}{2}\sigma} \widetilde{D (e^{-\frac{1}{2}\sigma} \cdot
\varphi)}, \] where $\,^{\sim} : S \longrightarrow \tilde{S}\,$
denotes the
canonical identification of the spinor bundles of $(M,g)$ and
$(M,\tilde{g})$.
\end{pr}
\vspace{0.1cm}

\begin{pr}  \label{pr3} {\em (\cite{Baum/Friedrich/ua:91} Cor.1.2)}\\
The dimension of the
space of twistor spinors is conformally invariant and bounded by \[
\dim \mbox{Ker} {\cal D} \le 2^{[\frac{n}{2}]+1} . \]
\end{pr}
\vspace{0.1cm}

\begin{pr} \label{pr4} {\em (\cite{Baum/Friedrich/ua:91} Cor.1.3)}\\
Let $\varphi \in \Gamma
(S)$ be a non-trivial twistor spinor and $x_0\in M$. Then $\,\varphi
(x_0)\neq 0\,$ or $\, D \varphi (x_0) \neq 0$.
\end{pr}
\ \\
Let $R$ be the scalar curvature and Ric the Ricci curvature of
$(M^{n,1},g)$. If $\,\dim M = n \ge 3,\,\,\,K\,$ denotes the (2,0)
-Schouten tensor
\[ K(X,Y) = \frac{1}{n-2} \left\{ \frac{R}{2(n-1)} g
- \mbox{ Ric} \right\}. \]
We always identify $TM$ with $TM^*$ using the metric
$g$. For a ($2,0$)-tensor field $B$ we denote by the same symbol $B$
the corresponding $(1,1)$-tensor field $B:TM \longrightarrow TM\,$,
$\;g(B(X),Y) = B(X,Y).\,$
Let $C$ be the (2,1)-Schouten-Weyl tensor \[
C(X,Y) = (\nabla_X K)(Y) - (\nabla_Y K)(X) . \]
Furthermore, let
$W$ be the (4,0)-Weyl tensor of $(M,g)$ and let denote by the same
symbol the corresponding (2,2)-tensor field $\; W : \Lambda^2 M
\longrightarrow \Lambda^2 M.\;$  Then we have

\begin{pr}\label{pr5} {\em (\cite{Baum/Friedrich/ua:91} Th.1.3, Th.1.5)}\\
Let $\varphi \in
\Gamma (S)$ be a twistor spinor and $\eta = Y \wedge Z \in
\Lambda^2M\,$ a two form. Then
\begin{eqnarray}
 D^2 \varphi & = &
\frac{1}{4} \frac{n}{n-1} R \varphi\,, \label{9} \\
\label{10} \nabla^S_X D
\varphi & = & \frac{n}{2} K(X) \cdot \varphi\;,\\
\label{11} W (\eta )
\cdot \varphi &=& 0 \;,\\
\label{12}
W(\eta)\cdot D\varphi& = & n\, C(Y,Z) \cdot \varphi \;,\\
\label{13}
(\nabla_X W) (\eta)\cdot\varphi & = & X\cdot C(Y,Z)\cdot\varphi
+\frac{2}{n}(X\ecke W(\eta))\cdot D\varphi\;. \hspace{4cm}
\end{eqnarray}
\end{pr}
\ \\
If the scalar curvature $R$ of $\,(M^{n,1},g)\,$ is constant and
non-zero, equation (\ref{9}) shows that the spinor fields
\[ \psi_{\pm} := \frac{1}{2} \varphi \pm
\sqrt{\frac{n-1}{nR}}\,D\varphi \]
are formal eigenspinors of the Dirac operator $D$ to the eigenvalue
$\,\pm \frac{1}{2} \sqrt{\frac{nR}{n-1}}\,$.\\
A special class of twistor spinors are the so-called {\em Killing
spinors} $\,\varphi \in \Gamma(S)\,$ defined by the condition
\[ \nabla_X^S \varphi = \lambda \,X \cdot \varphi \qquad \mbox{ for
all }\;\; X \in \Gamma(TM),\]
where $\lambda\,$ is a constant complex number, called the {\em Killing
number} of $\varphi$. Using the twistor equation and the properties
(9) and (10) one obtains that for an Einstein space $\,(M^{n,1},g)\,$
with constant scalar curvature $\,R \not = 0\,$ the spinor fields
$\,\psi_{\pm}\,$ are Killing spinors to the Killing number $ \,\lambda
= \mp \frac{1}{2} \sqrt{\frac{R}{n(n-1)}}\,$. Hence, on this class of
Lorentzian manifolds each twistor spinor is the sum of two Killing
spinors. Therefore, we are specially interested in non-Einsteinian
Lorentzian manifolds which admit twistor spinors.\\
To each spinor field we associate a vector field
in the following way.

\begin{de} Let $\varphi \in \Gamma (S)$. The vector filed
$V_{\varphi}$ definied by \[ g(V_{\varphi}, X):=-\langle X \cdot
\varphi, \varphi \rangle \,,\quad\qquad X \in \Gamma(TM) \]
is called the canonical vector field of $\varphi$.
\end{de}
Because of (1),
$V_{\varphi}$ is a real vector field. By Zero$( \varphi )$ and
Zero$(X)$
we denote the zero sets of a spinor field $\varphi$ or a vector
field $X$.

\begin{pr} \label{pr6}
\begin{enumerate}
\item  For each spinor
field $\varphi \in \Gamma (S)$ $\; \mbox{Zero}(\varphi ) = \mbox{
Zero}(V_{\varphi}). $
\item  If $n$ is even, $n \le 6$ and $\varphi
\in \Gamma (S^{\pm} ) $ is a half spinor, then $\; V_{\varphi} \cdot
\varphi = 0. \;$ In particular, $V_{\varphi}$ is an isotropic vector
field. \end{enumerate}
\end{pr}
{\bf Proof:} Let $\varphi \in \Gamma (S)$. From (5) follows for the
time orientation
$\xi$
\bean g(V_{\varphi},\xi)&=&
-\langle \xi \cdot \varphi , \varphi \rangle = - (\xi \cdot \xi \cdot
\varphi , \varphi )_{\xi} = -(\varphi, \varphi )_{\xi} . \eean
Since the scalar product $\,(\cdot,\cdot)_{\xi}\,$ is positive definite,
this shows that Zero$(V_{\varphi}) = $ Zero$(\varphi )$. The second
statement is proved by a direct calculation using a basis
representation
of $\varphi$ and $V_{\varphi}$ and the formulas (1) and
(2).\qed
In the Riemannian case Proposition 6.1 is not true.
There exist non-trivial spinor fields $\varphi$ such that the
canonical vector field $V_{\varphi}$ is identically zero (see
\cite{Kuehnel/Rademacher:95}). On the other hand, the zero set
Zero$(\varphi)$ of a
Riemannian twistor spinor is discret (\cite{Baum/Friedrich/ua:91},
Th.2.1). This is in the Lorentzian setting not the case.\\
We call a subset $A \subset M$ isotropic, if each
differentiable curve in $A$ is isotropic.

\begin{pr} \label{pr7} Let $\varphi \in \Gamma (S)$ be a twistor spinor.
Then the zero set of $\varphi$ is isotropic.
\end{pr}
{\bf Proof:} Let $\gamma : I \longrightarrow \mbox{Zero}(\varphi
)$ be a curve in Zero$(\varphi )$. Then $\,\varphi(\gamma(t))
\equiv 0\,$ and therefore $\nabla_{\dot{\gamma}(t)}\varphi\equiv 0$.
From the twistor equation (6) it follows $\,\dot{\gamma}(t)\cdot
D\varphi(\gamma(t)) \equiv 0\,$. Since by Proposition 4
$\,D\varphi(\gamma(t)) \not = 0\,$, $\,\dot{\gamma}(t)$ is isotropic
for all $t \in I$. \qed

\begin{pr} \label{pr8} Let $\varphi \in \Gamma (S)$ be a twistor spinor.
Then $V_{\varphi}$ is a conformal vector field and for the Lie
derivative  \[ L_{V_{\varphi}}
g = - \frac{4}{n} \,\mbox{Re} \langle \varphi , D \varphi \rangle\, g
 \]
holds.
\end{pr}
{\bf Proof:} Let $\,V:=V_{\varphi}\,$. From the definition of
$V_{\varphi}$ it follows
\bean (L_V g ) (X,Y) &=& g (
\nabla_X V, Y ) + g (X, \nabla_Y V ) \\
&=& X (g (V, Y )) - g (V, \nabla_X
Y ) + Y (g (X, V) )  - g ( \nabla_Y X, V ) \\
&=&- X \langle Y
\cdot \varphi, \varphi \rangle - Y \langle X \cdot \varphi, \varphi
\rangle  +\langle \nabla_X Y \cdot \varphi , \varphi \rangle +
\langle \nabla_Y X \cdot \varphi , \varphi \rangle \\ &
\stackrel{(\ref{4})}{=}
 &-\langle \nabla_X Y \cdot \varphi, \varphi
\rangle- \langle Y \cdot \nabla^S_X \varphi, \varphi \rangle -\langle
Y \cdot \varphi , \nabla^S_X \varphi \rangle -
\langle \nabla_Y X \cdot \varphi , \varphi \rangle \\
&&-\langle X
\cdot \nabla^S_Y \varphi, \varphi \rangle - \langle X \cdot \varphi ,
\nabla^S_Y \varphi \rangle + \langle \nabla_X Y \cdot \varphi ,
\varphi \rangle + \langle \nabla_Y X \cdot \varphi , \varphi \rangle
\\
& \stackrel{(\ref{3})}{=} &
- \langle Y \cdot \nabla^S_X \varphi + X \cdot \nabla^S_Y \varphi,
\varphi \rangle - \langle \varphi, Y \cdot \nabla^S_X \varphi + X
\cdot \nabla^S_Y \varphi \rangle\,. \eean
Using (\ref{7}) we obtain
\[ \left(
L_V g \right) (X,Y) = - \frac{4}{n} g (X,Y)\;\mbox{Re} \langle \varphi , D
\varphi \rangle . \] \qed
From Proposition 8 follows that for each twistor spinor $\varphi\;$
div$(V_{\varphi}) = - 2\, \mbox{Re} \langle \varphi, D \varphi
\rangle$. For the imaginary part of $\, \langle \varphi , D \varphi
\rangle\,$ we have

\begin{pr} \label{pr9} Let $\varphi \in \Gamma (S)$ be a twistor
spinor.
Then the function $\,C_{\varphi} := \,\mbox{Im}\, \langle \varphi , D
\varphi \rangle\,$ is constant on $M$.
\end{pr}
{\bf Proof:} Because of (\ref{3}) the function $\,\langle Y \cdot \psi, \psi
\rangle\,$ is real for each vector field $Y$ and each spinor field $\psi$.
Furthermore,
\bean X \langle D \varphi, \varphi \rangle & \stackrel{(\ref{4})}{=}
& \langle \nabla^S_X D
\varphi , \varphi \rangle + \langle D \varphi , \nabla^S_X \varphi
\rangle \\ & \stackrel{(\ref{6}),(\ref{10})}{=} &
 \frac{n}{2}\, \langle K(X)
\cdot\varphi,\varphi\rangle-\frac{1}{n}\,\langle D\varphi,X\cdot
D\varphi\rangle. \eean
Hence $\,X \langle D \varphi , \varphi \rangle\,$
is a real function. Therefore, $\,C_{\varphi} =\, \mbox{Im} \langle
\varphi , D \varphi \rangle$ is constant. \qed
Let us denote by $C$ the $(3,0)$-Schouten-Weyl tensor $\; C (X, Y, Z )
= g (X, C (Y,Z) ).$ \\

\begin{pr} \label{pr10} Let $\varphi \in \Gamma (S)$ be a twistor
spinor. Then
\begin{enumerate}
\item $V_{\varphi} \ecke C = 0 .$
\item If $n = 4$, then $\, V_{\varphi} \ecke W = 0. $
\end{enumerate}
\end{pr}
{\bf Proof:} From (\ref{11}) and (\ref{12}) we obtain
\bean C (V_{\varphi} , X, Y)
&=& g(V_{\varphi} , C ( X, Y )) = -\langle C ( X, Y) \cdot
\varphi , \varphi \rangle \\ &=& - \frac{1}{n}\, \langle W (X \wedge Y
) \cdot \varphi , \varphi \rangle = \frac{1}{n} \,\langle \varphi , W
(X \wedge Y ) \cdot \varphi \rangle \;=\; 0 . \eean
Let $\,\varphi = a
u (\varepsilon,1)+bu(-\var,-1)\in\Gamma (S^{\var})\,$ be a half spinor
on a 4-dimensional manifold. Then by a direct calculation using
(\ref{1}) and (\ref{2}) we obtain
\bean V_{\varphi} &=& (|a|^2 + |b|^2 ) s_1
+(|a|^2-|b|^2)s_2 - 2 \mbox{Re}
(ia\bar{b})s_3-2\var\mbox{Re}(a\bar{b})s_4. \eean
Hence,
\bea\label{12a} W(V_\varphi,s_i,s_j,s_k) &=&
(|a|^2+|b|^2)W_{1ijk}+(|a|^2-|b|^2)W_{2ijk}\nonumber\\ && -
2\mbox{Re}(ia\bar{b})W_{3ijk}-2\var\mbox{Re}(a\bar{b})W_{4ijk}. \eea
On the other hand, from the basis representation of
\[ 0 = W(s_j \wedge s_k)\cdot \varphi= \sum\limits_{r<l} \var_r \var_l
W_{rljk}\,
s_r \cdot s_l \cdot \varphi \]
result the equations
\begin{eqnarray}\label{13a}
0 & =& (W_{12jk}-\var iW_{34jk}) a + (i W_{13jk} - \var W_{24jk} -\var
W_{14jk} + i W_{14jk} ) \cdot b \\
\label{14a} 0 &=&
(-W_{12jk} + \var iW_{34jk}) b  + (-iW_{13jk}+\var W_{24jk} - \var
W_{14jk}+i W_{23jk}) a\,. \end{eqnarray}
Then looking at the real and
imaginary part of the equations $\,(\ref{13a}) \bar{a} \pm (\ref{14a}) \bar{b}\,$ and
$\,(\ref{13a}) \bar{b} \pm (\ref{14a})\bar{a}\,$ one obtains $\,W(
V_{\varphi}, s_i , s_j , s_k) = 0$.\qed

\section{Pseudo-hermitian geometry} \label{s4}

Before we define the Fefferman spaces we recall some basic facts from
pseudo-hermitian geometry in order to fix the notations. The proofs of
the following propositions are obtained by easy direct calculations
(see \cite{Tanaka:75}, \cite{Baum:97}).\\
Let $M^{2n+1}$ be a smooth connected manifold of odd dimension $2n+1$. A
{\em complex CR-structure} on $M$ is a complex subbundle $T_{10}$ of
$TM^{\Bbb C}$ such that\\[0.2cm]
\begin{tabular}{lll}
&1. & $ \dim_{\Bbb C} T_{10}=n,$\\[0.1cm]
&2. & $ T_{10}\cap\overline{T_{10}}=\{0\},$\\[0.1cm]
&3. & $ [\Gamma(T_{10}),\Gamma(T_{10})]\subset\Gamma(T_{10})\quad$
(integrability condition).
\end{tabular}\\[0.2cm]
A {\em real CR-structure} on $M$ is a pair $(H,J)$, where\\[0.2cm]
\begin{tabular}{lll}
& 1.& $ H\subset TM$ is a real $2n$-dimensional subbundle,\\
& 2.& $ J:H\longrightarrow H$ is an almost complex structure on $H:
\;J^2=-\mbox{id}$,\\[0.1cm]
& 3.& $\mbox{If }\,X,Y\in\Gamma(H)\,$, then
$\,[JX,Y]+[X,JY] \in\Gamma(H)\,$ and\\[0.1cm]
& & $ N_J(X,Y) := J([JX,Y]+[X,JY])-[JX,JY]+[X,Y] \equiv 0\,$\\
& &  (integrability condition).
\end{tabular}\\[0.2cm]
Obviously the complex and real CR-structure correspond to each other:
If $\,T_{10}\subset TM^{\Bbb C}\,$ is a complex CR-structure, then
$\,H:=\mbox{Re }(T_{10} \oplus\overline{T_{10}}))\,$,
$\,J(U+\bar{U}):= i(U-\bar{U})\,$ defines a real CR-structure. If
$(H,J)$ is a real CR-structure, then the eigenspace of the complex
extension of $J$ on $H^{\Bbb C}$ to the eigenvalue $i$ is a complex
CR-structure. A {\em CR-manifold} is an odd-dimensional
manifold equipped with a (real or complex) CR-structure. Let
$(M,T_{10})$ be a CR-manifold. The hermitian form on $T_{10}$
\[ L:T_{10}\times T_{10}\longrightarrow
E := TM^{\Bbb C}/_{T_{10}\oplus\overline{T_{10}}}\]
\[ L(U,V):= i[U,\bar{V}]_E\,, \]
where $X_E$ denotes the projection of
$X\in TM^{\Bbb C}$ onto $E$, is called the {\em Levi-form} of
$(M,T_{10})$. The CR-manifold is called {\em non-degenerate},
if its Levi-form $L$ is non-de\-ge\-ne\-ra\-te. An nowhere vanishing
1-form
$\theta \in\Omega^1(M)$ is called a {\em pseudo-hermitian structure}
on $(M,T_{10})$, if
$\,\theta|_H\equiv 0\,$.
$(M,T_{10},\theta)$ is called a {\em pseudo-hermitian manifold}.
There exists a pseudo-hermitian structure $\theta$ on $(M,T_{10})$ if
and only if $M$ is orientable. Two pseudo-hermitian structures
$\theta,\tilde{ \theta}$ differs by a real nowhere vanishing function
$\,f \in C^{\infty}(M)\,$: $\,\tilde{\theta}=f\cdot \theta\,$.
Let $(M,T_{10},\theta)$ be a pseudo-hermitian manifold. The hermitian
form $\;L_\theta:T_{10}\times T_{10}\longrightarrow \C\;$
\[ L_\theta(U,V):= -id\theta(U,\bar{V}) \]
is called the {\em Levi-form} of $(M,T_{10},\theta)$.
Obviously, we have $\,\theta(L(U,V))=L_\theta(U,V)\,$. The
pseudo-hermitian manifold $(M,T_{10},\theta)$ is called
{\em strictly pseudoconvex}, if the Levi-form $L_\theta$ is
positive definite. If the pseudo-hermitian manifold
$(M,T_{10},\theta)$ is non-degenerate, then the pseudo-hermitian
structure $\theta$ is a contact form. We denote by $T\in \Gamma
(TM)\,$ the {\em characteristic vector field} of this contact
form, e.g.
the vector field uniquely defined by
\[ \theta(T) \equiv 1 \qquad \mbox{ and } \qquad T
\;_-\!\rule{0.2mm}{0.2cm}\;\; d\theta \equiv  0. \]
From now on
we always suppose, that $(M,T_{10},\theta)$ is non-degenerate. If $M$
is oriented, we always choose $\theta$ such that a basis of the form
$\,(X_1,JX_1, \ldots , X_n, JX_n, T)\,$ is positive oriented on $M$.
We consider the following spaces of forms: \bean
\Lambda^{q,0}M &:=& \{\omega\in\Lambda^qM^{\Bbb C}\mid V\ecke \omega= 0\;
\;\;\forall V\in\overline{T_{10}}\}\\
\Lambda^{0,q}M &:=&
\{\omega\in\Lambda^qM^{\Bbb C}\mid V\ecke \omega=0\;\;\; \forall V\in
T_{10}\}\\
\Lambda^{p,q}M &:=&
\mbox{span}\{\omega\wedge\sigma\mid\omega\in\Lambda^{p,0} M,\;
\sigma\in\Lambda^{0,q}M\}\\
\Lambda^{p,q}_\theta M &:=& \{\omega\in\Lambda^{p,q}M\mid T\ecke
\omega=0\}.
\eean
Now, let us extend the Levi-form of $(M,T_{10},\theta)$ to $TM^{\Bbb C}$
by \\[0.2cm]
\hspace*{1cm} $ L_\theta(\bar{U},\bar{V}):=
\overline{L_\theta(U,V)}=L_\theta(V,U)\,, \quad L_\theta(U,\bar{V}):=0
\,,\quad$ where
$U,V\in T_{10}$,\\
\hspace*{1cm} $ L_\theta(T,\,\cdot\,):= 0. $

\begin{pr} \label{pr11} Let $L_\theta:TM^{\Bbb C}\times
TM^{\Bbb C} \longrightarrow \C$ be the Levi-form of $(M,T_{10},\theta)$
and let $T$ be the characteristic vector field of $\theta$. Then
\begin{eqnarray}
& & [T,Z]\in\Gamma(T_{10}\oplus
\overline{T_{10}})\qquad \;
\mbox{ if }\; Z\in\Gamma(T_{10})\, \mbox{ or }\, Z\in
\Gamma(\overline{T_{10}})\,, \label{14} \hspace{3.5cm} \\[0.1cm]
& &
L_\theta([T,U],V)+L_\theta(U,[T,V])=T(L_\theta(U,V)) \qquad
\forall \; U,V\in\Gamma (T_{10}) \,, \label{15}\\[0.1cm]
& &  L_\theta([T,\bar{U}],V)=L_\theta([T,\bar{V}],U)  \hspace{3.3cm}
\forall \; U,V\in
\Gamma(T_{10})\,,\label{16} \\[0.1cm]
& &  L_\theta([T,U],\bar{V})=L_\theta([T,V],\bar{U}) \hspace{3.3cm}
\forall \; U,V\in\Gamma (T_{10})\,,\label{17}
\end{eqnarray}
\end{pr}
\ \\
If we consider the Levi-form $L_\theta$ as a bilinear form on the real
tangent bundle, we obtain a symmetric bilinear form on $TM$ which is
non-degenerate on $H$.

\begin{pr} \label{pr12} Let $(M^{2n+1},T_{10},\theta)$ be
a non-degenerate pseudo-hermitian manifold and $(H,J)$ the real
CR-structure, defined by $T_{10}$. Let $X$ and $Y$ be two vector
fields in $H$. Then the Levi-form
$L_\theta:TM\times TM \longrightarrow \R$ satisfies
\begin{eqnarray}
& &
L_\theta(X,Y)= d\theta(X,JY)\,,\label{18}\\[0.1cm]
& & L_\theta(JX,JY)=L_\theta(X,Y) \quad \mbox{and} \quad
L_\theta(JX,Y)+L_\theta(X,JY)=0\,, \label{19} \hspace{3,5cm} \\[0.1cm]
& & L_\theta([T,X],Y)-L_\theta([T,Y],X) \,=\,
L_\theta([T,JX],JY)-L_\theta([T,JY],JX)\,.\label{20}
\end{eqnarray}
\end{pr}
\ \\
On non-degenerate pseudo-hermitian manifolds there exists a special
covariant derivative, the so-called {\em Webster connection}, which
was introduced by Tanaka (\cite{Tanaka:75} and by Webster
(\cite{Webster:78}).

\begin{pr} \label{pr13} Let $(M,T_{10},\theta)$ be a
non-degenerate pseudo-hermitian manifold and let $T$ be the
characteristic vector field of $\theta$. Then there exists an uniquely
determined covariant derivative $\;
\nabla^W:\Gamma(T_{10})\longrightarrow \Gamma(T^*M^{\Bbb C}\otimes
T_{10})\;$  on
$\,T_{10}\,$ such that\\[0.2cm]
1. $\;\nabla^W$ is metric with respect to $L_\theta:$
\vspace{-0.5cm}\\
\begin{eqnarray}
X(L_\theta(U,V))=L_\theta(\nabla^W_XU,V)+L_\theta(U,\nabla^W_{\bar{X}}
V) \qquad
U,V\in\Gamma(T_{10}), \; X\in \Gamma(TM^{\Bbb C}) \label{21} \end{eqnarray}
\vspace{-0.8cm}\\
\parbox{13cm}{
2.$ \quad \nabla^W_TU=\mbox{pr}_{10} [T,U]$, \\[0.2cm]
3.$\quad \nabla^W_{\bar{V}}U =\mbox{pr}_{10}[\bar{V},U]$,} \hfill
\parbox{6mm}{\begin{eqnarray} \label{22}\\[0.1cm] \label{23}
\end{eqnarray}}
\vspace{-0.2cm}\\ where
$\mbox{pr}_{10}$ denotes the projection on $T_{10}\;$. Furthermore,
$\nabla^W$ satisfies \be\label{24}
\nabla^W_UV-\nabla^W_VU=[U,V],\qquad U,V\in\Gamma(T_{10}). \ee
\end{pr}
\ \\
Now, we extend the Webster connection to $TM^{\Bbb C}$ by
\bean
\nabla^W\bar{U} := \overline{\nabla^WU} \qquad \mbox{ and } \qquad
\nabla^WT := 0. \eean

\begin{pr} \label{pr14} The torsion $\mbox{Tor}^W$ of the Webster
connection $\; \nabla^W: \Gamma(TM^{\Bbb C})\longrightarrow
\Gamma(T^*M^{\Bbb C}\otimes
TM^{\Bbb C})\;$ satisfies
\begin{eqnarray}  \label{26}
\mbox{Tor}^W(U,V) &=& \mbox{Tor}^W(\bar{U},\bar{V})\,=\,0\,,
\hspace{5cm}\\
\label{27}
\mbox{Tor}^W(U,\bar{V}) &= &i L_\theta(U,V) \,T \,,\\
\label{28}
\mbox{Tor}^W(T,U) &=& -\mbox{pr}_{01}[T,U]\,,\\
\mbox{Tor}^W(T,\bar{U})
&=& -\mbox{pr}_{10}[T,\bar{U}]\,, \end{eqnarray}
where $\mbox{pr}_{01}$ denotes the
projection onto $\overline{T_{10}}\,$, $\,p_{10}$ the projection onto
$T_{10}\,$ and $U,V \in \Gamma(T_{10})$.
\end{pr}
\ \\
Let $(M,T_{10},\theta)$ be a non-degenerate pseudo-hermitian
manifold and let $(p,q)$ be the signature of $(T_{10},L_\theta)$. Then
$\; g_\theta:= L_\theta+\theta\circ \theta \;$ defines a metric of
signature $(2p,2q+1)$ on $M$.

\begin{pr} \label{pr15} Let $(M,T_{10},\theta)$ be a non-degenerate
pseudo-hermitian manifold. Then the Webster connection
$\nabla^W:\Gamma(TM) \longrightarrow \Gamma(T^*M\otimes TM)\,$
considered on the real
tangent bundle is metric with respect to $\,g_\theta\,$ and the
torsion of
$\nabla^W$ is given by
\bea
\label{29} \mbox{Tor}^W(X,Y)& = &
L_\theta(JX,Y)\cdot T \qquad \qquad \mbox{ for }\; X,Y\in\Gamma(H),
\hspace{2cm}\\
\mbox{Tor}^W(T,X) &= &
-\frac{1}{2}\{[T,X]+J[T,JX]\} \qquad \mbox{ for } \; X\in\Gamma(H).\label{30}
\eea
Furthermore, on $\Gamma(H)$ \be\label{31} \nabla^W\circ J=J\circ
\nabla^W\,. \ee
\end{pr}
\ \\
Now, let $\;R^{\nabla^W}\in\Gamma
(\Lambda^2M^{\Bbb C}\otimes\mbox{End}(TM^{\Bbb C}, TM^{\Bbb C}))\;$ be the curvature
operator of $\nabla^W$
\[ R^{\nabla^W}(X,Y) =
[\nabla^W_X,\nabla^W_Y]-\nabla^W_{[X,Y]}. \]
Then the (4,0)-curvature tensor ${\cal R}^W$
\[ {\cal R}^W(X,Y,Z,V):=
g_\theta(R^{\nabla^W}(X,Y)Z,\bar{W}),\qquad X,Y,Z,W\in TM^{\Bbb C} \]
has the following symmetry properties \\

\begin{pr} \label{pr16} Let $\,X,Y,Z,V\in TM^{\Bbb C}\,$, $\,A,B,C,D\in T_{10}\,$.
Then\\[0.2cm]
\hspace*{1cm} ${\cal R}^W(X,Y,Z,V)=-{\cal R}^W(Y,X,Z,V)=-{\cal
R}^W(X,Y,V,Z)$\\[0.1cm]
\hspace*{1cm}
$\overline{{\cal R}^W(X,Y,Z,V)}={\cal R}^W(\bar{X},\bar{Y},\bar{Z},
\bar{V})$\\[0.1cm]
\hspace*{1cm}
${\cal R}^W(A,\bar{B},C,\bar{D})={\cal
R}^W(C,\bar{B},A,\bar{D})$\\[0.1cm]
\hspace*{1cm}
${\cal R}^W(A,B, \cdot , \cdot )=0$
\end{pr}
\ \\
Let $\omega\in\Lambda^2M^{\Bbb C}$ be a complex 2-form and
$\,\tilde{\omega}: T_{10} \longrightarrow T_{10}\,$ the uniquely
determined $\C$-linear
map with $\,\omega(U,\bar{V})=L_\theta(\tilde{\omega}U,V)\,$,
$\,U,V \in T_{10}\,$. Then the $\theta$-trace of $\omega$ is
defined by $\;
\mbox{Tr}_\theta\omega:=\mbox{Tr}(\tilde{\omega}).\;$ If
$\,(Z_1,\ldots,Z_n)\,$ is an unitary basis of
$\,(T_{10},L_\theta)\,$,
$\,\varepsilon_k= L_\theta(Z_k,Z_k)\,$, then \[
\mbox{Tr}_\theta\omega=\suml^n_{\alpha
=1}\varepsilon_\alpha\;\omega(Z_\alpha,
\bar{Z}_\alpha). \]
The (2,0)-tensor field
\[ \mbox{Ric}^W:= \mbox{Tr}_\theta^{(3,4)}{\cal
R}^W=\suml^n_{\alpha=1}\var_\alpha {\cal R}^W( \cdot , \cdot ,
Z_\alpha,\bar{Z}_\alpha) \]
is called the
{\em Webster-Ricci-tensor}, the function
$\,\, R^W:=\mbox{Tr}_\theta\mbox{Ric}^W\,\,$
is the {\em Webster scalar curvature}. Proposition 16 shows that
$\,\mbox{Ric}^W\in\Lambda^{1,1}M\,$,
$\,\mbox{Ric}^W(X,Y) \in i\R\,$ for all $X,Y\in TM\,$ and
that $R^W$ is a real function.

\section{Fefferman spaces} \label{s5}

Let $\,(M^{2n+1},T_{10})\,$ be a CR-manifold. The complex line bundle
$\,K:=\Lambda^{n+1,0}M\,$ of $(n+1,0)$-forms is called the
{\em canonical bundle} of $\,(M^{2n+1},T_{10})\,$. $\R^+$ acts on
$\,K^*=K\setminus \{0\}\,$ by multiplication. Let
$\,F:=K^*/_{\R^+}\,$. Then
$(F,\pi,M)$ is the $S^1$-principal bundle over $M$ associated to $K$.
We call $(F,\pi,M)$ the {\em canonical $S^1$-bundle} of
$(M,T_{10})$.
Now, let $(M,T_{10},\theta)$ be a non-degenerate pseudo-hermitian
manifold and $\,\nabla^W:\Gamma(T_{10})\longrightarrow
\Gamma(T^*M^{\Bbb C}\otimes
T_{10})\,$ its Webster-connection. $\nabla^W$ allows us to define a
connection $A^W$ on the canonical $S^1$-bundle $F$ in the following
way: Let $\,s=(Z_1, \ldots ,Z_n)\,$
be a local unitary basis of $(T_{10},L_\theta)$ over $U \subset M$ and let
us denote by $\,\omega_s:=(\omega_{\alpha\beta})\,$ the matrix of
connection forms of $\nabla^W$ with respect to $s$
\[
\nabla^WZ_\alpha=\suml_\beta\omega_{\alpha\beta}Z_\beta. \]
\vspace{-0.3cm}\\
$(Z_1,\ldots,Z_n,\bar{Z}_1,\ldots,\bar{Z}_n,T)\,$ is a local basis of
$TM^{\Bbb C}$ over $U$. Let
$(\theta^1,\ldots,\theta^n,\bar{\theta}^1,\ldots,
\bar{\theta}^n,\theta)$ be the corresponding dual basis. Then
\[
\hat{\tau}_s := \theta\wedge\theta^1\wedge\ldots\wedge\theta^n : U
\longrightarrow K \]
is a local section in $K$. We denote by $\,\tau_s := [\hat{\tau}_s]\,$
the corresponding local section in $\,F= K^*/_{\R^+}\,$.
The Webster connection $\,\nabla^W\,$ defines in the standard way a
covariant derivative $\,\nabla^K\,$ in the canonical line bundle $K$
such that
\[ \nabla^K \hat{\tau}_s = - \sum\limits_{\alpha} \omega_{\alpha\alpha}
\cdot \hat{\tau}_s = - \mbox{ Tr}\,\omega_s \cdot \hat{\tau}_s \,.\]
\vspace{-0.3cm}\\
Since $\,\nabla^W\,$ is metric with respect to $L_\theta$,
the trace Tr$\,\omega_s\,$ is purely imaginary. Hence $\nabla^K\,$ is induced
by a connection $A^W\,$ on the associated $S^1$-principal bundle
$\,(F,\pi,M;S^1)\,$ with the local connection forms
\[ \tau_s^* A^W = - \mbox{ Tr}\, \omega_s \,.\]
Let $\Omega^W$ be the curvature form of the
connection $A^W$ on $F$. Since $\Omega^W$ is tensionell and
right-invariant, it can be
considered as 2-form on $M$ with values in $i\R$.
Over $U \subset M$
\be\label{33} \Omega^W=dA^{\tau_s}=-\mbox{Tr}\,d\omega_s. \ee
holds. On
the other hand,
\bean \mbox{Ric}^W(X,Y) &=&\suml_\alpha\var_\alpha
L_\theta(([\nabla^W_X,\nabla^W_Y]
-\nabla^W_{[X,Y]})Z_\alpha,\bar{Z}_\alpha)\\ &=&\Big(\suml_\alpha
d\omega_{\alpha\alpha}-\suml_{\alpha,\beta}\omega_{\alpha
\beta}\wedge\omega_{\beta\alpha}\Big)(X,Y). \eean
\vspace{-0.5cm}\\
Hence, \bean
\mbox{Ric}^W =
\mbox{Tr}\,d\omega_s-\mbox{Tr}\,(\omega_s\wedge\omega_s) =
\mbox{Tr}\,d\omega_s. \eean
From (\ref{33}) it follows
\be \label{32} \Omega^W=-\mbox{Ric}^W\,. \ee
The connection $A^W$ on the canonical $S^1$-bundle $(F,\pi,M)$ is
called the {\em Webster-con\-nec\-tion on $F$}.
Two connections on an $S^1$-principal bundle over $M$ differ by an
1-form on $M$ with values in $i\R$. The connection
\[ A_\theta :=
A^W-\frac{i}{2(n+1)} R^W \theta \]
on the canonical $S^1$-bundle
$(F,\pi,M)$ is called the {\em Fefferman connection on $F$}.\\
Let us consider the following right-invariant metric on $F$: \[
h_\theta:=\pi^*L_\theta-i\frac{4}{n+2}\pi^*\theta\circ A_\theta, \]
where $\circ$ denotes the symmetric tensor product.
$h_\theta$ is
called the {\em Fefferman metric on $F$}. If $(T_{10},
L_\theta)$ is of signature $(p,q)$, then $h_\theta$ has signature
$(2p+1,2q+1)$. In particular, if $(M,T_{10},\theta)$ is strictly
pseudoconvex, $h_\theta$ is a Lorentzian metric.
The semi-Riemannian
manifold $(F,h_\theta)$ is called the {\em Fefferman space of
$\,(M,T_{10},\theta)\,$}. The fibres of the ca\-no\-nical $S^1$-bundle $F$ are
isotropic submanifolds of $(F,h_\theta)$. From the special choise of
the Fefferman connection $A_\theta$ in the definition of $h_\theta$
results
that the conformal class $[h_\theta]$ of the metric $h_\theta$ is an
invariant of the oriented CR-manifold $(M,T_{10})$, e.g. if
$\tilde{\theta}=f\cdot\theta$, $f>0$, is a further pseudo-hermitian
structure on $(M,T_{10})$, then $\,h_{\tilde{\theta}}=f\cdot
h_\theta\, $
(see \cite{Lee:86}, Th. 5.17.). We remark that Fefferman spaces are
never Einsteinian. \\
\ \\
In the following we always assume that
$(M,T_{10},\theta)$ is strictly pseudoconvex.
In order to find global solutions of the Lorentzian twistor equation
on Fefferman spaces it is necessary to change the topological type of
the canonical $S^1$-bundle.\\

\begin{pr} \label{pr17}
Let $(M^{2n+1},T_{10},\theta)$ be a strictly
pseudoconvex spin manifold. Then each spinor structure of the Riemann
manifold $(M,g_\theta)$ defines a square root $\sqrt{F}$ of the
canonical $S^1$-bundle $F$. (e.g. $\sqrt{F}$ is an $S^1$-bundle over
$M$ such that the associated line bundle $\,L:=\sqrt{F}
\times_{S^1}\C\,$ satisfies $\,L\otimes L=K\,$).
\end{pr}
{\bf Proof:} Let $\, U(n) \hookrightarrow SO(2n) \hookrightarrow
SO(2n+1)\,$ be the canonical embedding of $U(n)$ in $SO(2n+1)$.
\[
P_H:=\{(X_1,JX_1,\ldots,X_n,JX_n,T)\mid (X_1,JX_1,\ldots,X_n,JX_n)\;
\mbox{on-basis of } (H,L_\theta)\} \]
is an $U(n)$-reduction of the bundle $P_M$ of $SO(2n+1)$-frames of
$(M,g_\theta)$ .
Let $(Q_M,f_M)$ be a spinor structure of $(M,g_\theta)$ and let us
denote by $(Q_H,f_H)$ the reduced spinor structure \[ Q_H:=
f^{-1}_M(P_H),\qquad f_H:= f_M\mid_{Q_H}. \] Now, the proof of
Proposition \ref{pr17} is a repetition of Hitchin's proof of the fact that
each spinor structure on a K\"ahler manifold defines a square root of
the canonical bundle (see \cite{Hitchin:74}). Since we need some notation
later on, we repeat the idea of the proof.\\
Let $\ell:U(n) \longrightarrow
\mbox{Spin}(2n)^{\Bbb C}=\mbox{Spin}(2n) \times_{\Bbb Z_2}S^1$ be
defined by
\be\label{34}
\ell(A)=\prod\limits^n_{k=1}
\left(\cos\frac{\theta_k}{2}+\sin\frac{\theta_k}{2}
\cdot f_k\cdot J_0(f_k)\right) \times
e^{\frac{i}{2}\suml^n_{k=1}\theta_k}\,,
\ee
where $(f_1,\ldots,f_n)$ is an unitary basis of $\C^n$ such that
$Af_k= e^{i\theta_k}f_k$ and $J_0:\C^n \to \C^n$ is the standard
complex structure of $\C^n$. Then we have the following commutative
diagram
\vspace{-0.5cm}\\
\begin{center}
\setlength{\unitlength}{1cm}
\begin{picture}(10,3.3)
\put(0,0.3){ $S^1$ } \put(3,0.4){\vector(-3,0){1.5}} \put(2,0.5){det}
\put(4,0.3){ $U(n)$}
\put(6,0.4){\vector(3,0){1.5}} \put(6.7,0.5){$i$} \put(8,0.3){SO($2n$)}

\put(0,2.8){ $S^1$ } \put(1.5,2.9){\vector(3,0){1.5}}
\put(2,2.6){$j_2$}
\put(3.5,2.8){Spin$(2n)^{\Bbb C}$} \put(7.5,2.9){\vector(-3,0){1.5}}
\put(6.7,2.6){$j_1$} \put(7.9,2.8){Spin($2n$)}

\put(0.3,2.3){\vector(0,-2){1.3}} \put(4.5,1){\vector(0,2){1.3}}
\put(8.7,2.3){\vector(0,-2){1.3}}
\put(4.2,1.5){$\ell$} \put(8.9,1.5){$\lambda$}
\put(0.5,2){$z$}\put(0.5,1.6){$\downarrow$}\put(0.5,1.2){$z^2$}

\put(5.7,2.3){\vector(3,-2){2}} \put(5.8,1.4){{\footnotesize $\lambda
\circ pr_1$}} \put(0.54,1.8){-}

\end{picture}
\end{center}
where $i,j_1,j_2$ denote the canonical
embeddings and $\lambda:\mbox{Spin} (2n)\to SO(2n)$ is the universal
covering of $SO(2n)$. Hence, for each $A\in U(n)$ and each square
root of $\det(A)$ one has
\[ \lambda^{-1}(A):=
j_1\lambda^{-1}(i(A))=\pm\ell(A)\,\mbox{Det}(A)^{-\frac{1}{2}}. \]
Now,
let $\{(U_{\alpha\beta}\,,\, g_{\alpha\beta}:U_{\alpha\beta}\to
\lambda^{-1}(U(n)))\}_{\alpha,\beta}$ are the cocycles defining the
reduced spinor structure $\,(Q_H,f_H)\,$. Then on $U_{\alpha\beta}$ we
choose a square root $\;
h_{\alpha\beta}:U_{\alpha\beta}\to S^1 \;$
of the determinant of $\,\lambda(g_{\alpha\beta})^{-1}\,$ such that
\be\label{35}
h^2_{\alpha\beta}=\mbox{Det}(\lambda(g_{\alpha\beta}))^{-1}  \qquad
\mbox{ and } \qquad
 g_{\alpha\beta}=\ell(\lambda(g_{\alpha\beta}))\cdot
h_{\alpha\beta}. \ee
$\{(U_{\alpha\beta},h_{\alpha\beta})\}_{\alpha\beta}$ are cocyles
defining a square root $\,(\sqrt{F},\pi,M)\,$ of the canonical
$S^1$-bundle $(F,\pi,M)$. \qed
\ \\
Let $\,(\sqrt{F},\pi,M)\,$ be the square root of the canonical
$S^1$-bundle defined by the spinor structure of $(M,g_\theta)$. Then
the Webster connection $A^W$ on $F$ defines a corresponding connection
$A^{\sqrt{W}}$ on $\sqrt{F}$: Let $\,\{\tilde{s}_\alpha:U_\alpha\to
Q_H\}\,$ be a covering of $Q_H$ by local sections with the transition
functions $\,g_{\alpha\beta}\;$; $\;\tilde{s}_\alpha
=\tilde{s}_\beta\cdot
g_{\alpha\beta}\;$. Let $\,s_\alpha=f_H(\tilde{s}_\alpha) \in P_H\,$ and
denote by $\,\sqrt{\tau_{s_\alpha}}:U_\alpha\to\sqrt{F}\,$ the local
sections in $\sqrt{F}$ with transition functions $\,h_{\alpha\beta}\,$
\[
\sqrt{\tau_{s_\alpha}}=\sqrt{\tau_{s_\beta}}\cdot h_{\alpha\beta}, \]
defined by (\ref{35}). Then the local connection
forms of $A^{\sqrt{W}}$ are given by
\be\label{37}
\sqrt{\tau_{s_\alpha}}^{\,*}
A^{\sqrt{W}}=\frac{1}{2}\tau^*_{s_\alpha}A^W=-
\frac{1}{2}\,\mbox{Tr}\,\omega_{s_\alpha}
\ee
and the curvature of $A^{\sqrt{W}}$ is
\be\label{38}
\Omega^{\sqrt{W}}=\frac{1}{2}\Omega^W= - \frac{1}{2}\mbox{Ric}^W. \ee
The
connection $\,A^{\sqrt{}}_\theta\,$ on $\,\sqrt{F}\,$ defined by
\[ A^{\sqrt{}}_\theta:=
A^{\sqrt{W}}-\frac{i}{4(n+1)}R^W\cdot\theta \]
is called the
{\em Fefferman connection on $\sqrt{F}$} and the Lorentzian
metric \[ h_\theta:=\pi^*L_\theta- i\frac{8}{n+2}\pi^*\theta\circ
A^{\sqrt{}}_\theta \]
is the {\em Fefferman metric on $\sqrt{F}$}.
As we will see in the next section, the spinor structure $(Q_M,f_M)$
of $\,(M,g_\theta)\,$ defines a canonical spinor structure on
$\,(\sqrt{F},h_\theta)\,$.

\begin{de} The Lorentzian spin manifold $(\sqrt{F},h_\theta)$
is called the Fefferman space of the strictly
pseudoconvex spin manifold $\,(M,T_{10},\theta,(Q_M,f_M))\,$.
\end{de}
\ \\
\section{Spinor calculus for $S^1$-bundles with isotropic fibre over
strictly pseudoconvex spin manifolds} \label{s6}

Let $(M^{2n+1},T_{10},\theta)$ be a strictly pseudoconvex manifold and
let $(Q_M,f_M)$ be a spinor structure of $(M,g_\theta)$. Furthermore,
consider an $S^1$-principle bundle $(B,\pi,M;S^1)$ over $M$, a
connection $A$ on $B$ and a constant $c\in\R\backslash\{0\}$. Then \[
h:= h_{A,c} :=\pi^*L_\theta-i c\,\pi^*\theta\circ A \] is a Lorentzian
metric on $B$. In this section we want to derive a suitable spinor
calculus for the Lorentzian manifold $(B,h)$.\\
Let $\,N\in \Gamma(TB)\,$ be the fundamental vector field on $B$ defined
by the element $\frac{2}{c}i\in i\R$ of the Lie algebra $i\R$ of $S^1$
\[ N(b)=\widetilde{\frac{2}{c}i}\,(b):=\frac{d}{dt}\left(b\cdot
e^{\frac{2}{c}it}\right)|_{ t=0}. \]
Denote by $\,T^*\in \Gamma(TB)\,$ the
$A$-horizontal lift of the characteristic vector field $T$ of
$\theta$. Then $N$ and $T^*$ are global isotropic vector fields on $B$
such that $h(N,T^*)=1$. Consider the global vector fields
\be\label{39} s_1 = \frac{1}{\sqrt{2}}(N - T^*) \qquad \mbox{ and }
\qquad s_2 =
\frac{1}{\sqrt{2}} (N+T^*). \ee
Then
\[ h(s_1,s_1)=-1, \qquad
h(s_2,s_2)=1, \qquad h(s_1,s_2)=0. \]
Let the time orientation of
$(B,h)$ be given by $s_1$ and the space orientation by the vectors
$\,(s_2,X^*_1,JX_1^*,\ldots,X^*_n,JX_n^*))\,$, where
$\,(X_1,JX_1,\ldots,
X_n,JX_n)\in P_H\,,$ and $X^*$ denotes the $A$-horizontal lift of a vector
field $X$ on $M$. Now, let $(Q_H,f_H)$ be the reduced spinor structure
of $(M,g_\theta)$ defined in the previous section. Denote by
\[ S_H:=
Q_H \times_{\lambda^{-1}(U(n))}\Delta_{2n,0} \]
the corresponding spinor
bundle of $(H,L_\theta)$. Obviously, the bundle \[
\hat{P}_B:=\{(s_1,s_2,X^*_1,JX_1^*,\ldots,X^*_n,JX_n^*) \mid
(X_1,JX_1,\ldots,X_n,
JX_n) \mbox{ on-basis of } (H,L_\theta)\} \]
is an $U(n)$-reduction of
the frame bundle $P_B$ of $(B,h)$ with respect to the embedding
$\,U(n) \hookrightarrow SO_0(2n+2,1)\,$.
Since $\,\hat{P}_B\approx\pi^*P_H\,$ we have
\[
P_B \,\approx\, \pi^*P_H \times_{U(n)}\,SO_0(2n+2,1). \]
Therefore,
\bean Q_B
:= \pi^*Q_H \times_{\lambda^{-1}(U(n))}\,\mbox{Spin}_0(2n+2,1)\,, \qquad
f_B := [f_H,\lambda] \eean is a spinor structure of the Lorentzian
manifold
$(B,h)$. The corresponding spinor bundle $S$ on $(B,h)$ is given by
\be\label{41} S \,=\, \pi^*Q_H
\times_{\lambda^{-1}(U(n))}\,\Delta_{2n+2,1}.
\ee
\vspace{0.2cm}
\begin{pr} \label{pr19} Let $S_H$ be the spinor bundle of $(H,L_\theta)$
over $M$. Then the spinor bundle $S$ of $(B,h)$ can be identified with
the sum
\[ S \approx \pi^*S_H\oplus \pi^*S_H, \]
where the Clifford multiplication is given by
\begin{eqnarray}
\label{42} s_1\cdot
(\varphi,\psi) &=& (-\psi,-\varphi) \\
\label{43} s_2\cdot (\varphi,\psi) & = & (-\psi,\varphi)  \\
\label{44}
X^*\cdot(\varphi,\psi) &=& (-X\cdot\varphi,X\cdot\psi),\qquad X\in H.
\end{eqnarray}
In particular,
\begin{eqnarray}
\label{45} N\cdot(\varphi,\psi) & =& (-\sqrt{2}\,\psi,0)  \\
\label{46} T^*\cdot(\varphi,\psi) & = & (0,\sqrt{2}\,\varphi).
\end{eqnarray}
Furthermore, the positive and negative parts of $S$ are
\be\label{47} S^+ = \pi^*S^+_H\oplus\pi^*S^-_H,\qquad
S^- = \pi^*S^-_H \oplus \pi^*S^+_H.
\ee
The indefinite scalar product
$\,\langle \cdot , \cdot \rangle\,$ in $S$ is given by
\be\label{48}
\langle(\varphi,\psi),
(\hat{\varphi},\hat{\psi})\rangle=-(\psi,\hat{\varphi})_{S_H}
-(\varphi,\hat{\psi})_{S_H},
\ee
where $(\cdot,\cdot)_{S_H}$ is the
usual positive definite scalar product in $S_H$.
\end{pr}
{\bf Proof:} By definition of the spinor bundle $S$  (see (\ref{41}))
we have only to check, how the $\,\mbox{Spin}(2n)$-modul
$\Delta_{2n+2,1}$ decomposes into
$\mbox{Spin}(2n)$-representations. Let the embedding
$i:\R^{2n}\to\R^{2n+2,1}$ be given by $i(x)=(0,0,x)$ and let
$\mbox{Spin}(2n)\hookrightarrow \mbox{Spin}_0(2n+2,1)$ be the
corresponding embedding of the spin groups. Consider the following
isomorphisms of the representation spaces
\begin{eqnarray*}
\begin{array}{lccrcl}
 \chi :& \Delta_{2n+2,1} & \longrightarrow &
\,\Delta_{2n,0} &\oplus & \Delta_{2n,0}\, \\[0.2cm]
& \, u \otimes u(1) + v \otimes
u(-1)\, &\longmapsto &(u&,&v)
\end{array}
\end{eqnarray*}
where we use the notation of section \ref{s2}. Then
formula (\ref{1}) shows that
\begin{eqnarray*}
\begin{array}{lcl}
 \chi\,(e_1\cdot (u\otimes u(1)+v\otimes u(-1)))& =&
(-u,-v) \\[0.1cm]
\chi\,(e_2\cdot(u\otimes u(1) + v\otimes u(-1)))&=&(u,-v)\\[0.1cm]
\chi\,(e_k\cdot(u\otimes u(1) + v\otimes u(-1)))&=&(-e_{k-2}\cdot
u,e_{k-2}\cdot v), \qquad k>2.
\end{array}
\end{eqnarray*}
Therefore, $\chi$ is an isomorphism of the
$\mbox{Spin}(2n)$-representations and (\ref{42})-(\ref{44}) and
because of (\ref{39}) also the formulas (\ref{45}),
(\ref{46}) are valid. Let $\omega_{2n+2}=e_1\cdot\ldots\cdot
e_{2n+2}$ be the volume element of $\mbox{Cliff}_{2n+2,1}$ and
$\omega_{2n}=e_1\cdots e_{2n}$ the volume element of
$\mbox{Cliff}_{2n,0}$. Then using the identification $\chi$ we obtain
\[
\omega_{2n+2} \cdot (u,v) = (-\omega_{2n} \cdot u\,,\,\omega_{2n}\cdot
v). \]
According to the definition of $S^\pm$ this shows (\ref{47}). Because
of (\ref{5}) the scalar product satisfies
\bean
\langle(\varphi,\psi),(\hat{\varphi},\hat{\psi})\rangle &=& (s_1\cdot
( \varphi,\psi), (\hat{\varphi},\hat{\psi}))_{s_1} \\&=&
((-\psi,-\varphi), (\hat{\varphi},\hat{\psi}))_{s_1}\\ &=&
-(\psi,\hat{\varphi})_{S_H}-(\varphi,\hat{\psi})_{S_H}.
\eean
\qed
\ \\
In order to describe the spinor derivative in the spinor bundle $S$ of
$B$ we need the
connection forms of the Levi-Civita connection of $(B,h)$.  Let
$X,Y,Z$ be local vector fields on $(B,h)$ of constant length and
constant scalar products with each other. Then the Levi-Civita connection
$\nabla$ of $(B,h)$ satisfies
\be\label{49}
h(\nabla_XY,Z)=\frac{1}{2}\{h([X,Y],Z)+h([Z,Y],X)+h([Z,X],Y)\}.
\ee
For a vector $Z\in T_bB$ we denote by $Z^h$ the projection on the
horizontal tangent space and by $Z^v$ the projection on the vertical
tangent space. If $X\in T_{\pi(b)}M$, then $X^*\in T_bB$ denotes the
horizontal lift of $X$. Let $\,\Omega^A \in \Omega^2(M;i\R)\,$ be the
curvature form of the connection $A$. From the connection theory in
principle bundles follows for vector fields $X,Y$ on $M$
\begin{eqnarray}
\label{50} [X^*,N]\,\,\, & = & 0 \,,\\
\label{51}
[X^*,Y^*]^v &= & i\,\frac{c}{2}\,\Omega^A(X,Y)\cdot N\,,\\
\label{52}  [X^*,Y^*]^h &= & [X,Y]^*\,.
\end{eqnarray}
Now, let $X,Y\in\Gamma(H)$. Since $\;[T,X] \in
\Gamma(H)\;$ and \[
[X,Y] = \mbox{pr}_H[X,Y] + \theta([X,Y])\cdot T
= \mbox{pr}_H[X,Y]-d\theta(X,Y)\cdot T
\]
we obtain from
(\ref{51}) and (\ref{52})
\bea
\label{53} [T^*,X^*] &=&
[T,X]^* + i\,\frac{c}{2}\,\Omega^A(T,X) \cdot N\,, \\
\mbox{}[X^*,Y^*]
&=& \mbox{pr}_H [X,Y]^* - d\theta (X,Y) \cdot T^* +
i\,\frac{c}{2}\,\Omega^A(X,Y) \cdot N. \label{54}
\eea
\\
\begin{pr} \label{pr20}
Let $X,Y,Z\in\Gamma(H)$ be vector fields of constant lenght and
constant $L_\theta$-scalar products with each other. Then
\bean
h(\nabla_{X^*}Y^*,Z^*) &=& L_\theta(\nabla^W_XY,Z)\\
h(\nabla_NY^*,Z^*) &=& \frac{1}{2} d\theta(Y,Z)\\
h(\nabla_{T^*}Y^*,Z^*) &=&
\frac{1}{2}\{L_\theta([T,Y],Z)-L_\theta([T,Z],Y)
-i\frac{c}{2}\Omega^A(Y,Z)\}\\
h(\nabla_{X^*}Y^*,N) &=& -\frac{1}{2}d\theta(X,Y)\\
h(\nabla_{X^*}Y^*,T^*) &=&
\frac{1}{2}\{L_\theta([T,X],Y)+L_\theta([T,Y],X) +
i\frac{c}{2}\Omega^A(X,Y)\}\\
h(\nabla_{T^*}T^*,Z^*) &=& -i\,\frac{c}{2}\,\Omega^A(T,Z)\\
h(\nabla N,T^*) &=& h(\nabla
T^*,T^*)\,=\,h(\nabla N^*,N^*)\,=\,0 \\
h(\nabla_NN,Z^*) &=&
h(\nabla_NT^*,Z^*)\,=\,h(\nabla_{T^*}N,Z^*)\,=\,0.
\eean
\end{pr}
{\bf Proof:} From (\ref{49}) and (\ref{54}) it follows
\bean
2\,h(\nabla_{X^*}Y^*,Z^*) &=&
h(\,\mbox{pr}_H[X,Y]^*,Z^*)+h(\,\mbox{pr}_H[Z,Y]^*, X^*) +
h(\,\mbox{pr}_H[Z,X]^*,Y^*)\\ &=&
L_\theta([X,Y],Z)+L_\theta([Z,Y],X)+L_\theta([Z,X],Y).
\eean
According to (\ref{29}) $\,\mbox{Tor}^W(X,Y)=L_\theta(JX,Y)\cdot T\,$.
Hence,
\bean
L_\theta([X,Y],Z) &=&
L_\theta(\nabla^W_XY-\nabla^W_YX-\mbox{Tor}^W(X,Y),Z)\\ &=&
L_\theta(\nabla^W_XY-\nabla^W_YX,Z).
\eean
Therefore, using that
$\nabla^W$ is metric with respect to $L_\theta$ we obtain
\[
h(\nabla_{X^*}Y^*,Z^*)=L_\theta(\nabla_X^WY,Z).
\]
The other formulas
follow immediately from the definition of $h$ and (\ref{49}),
(\ref{50}), (\ref{53}) and (\ref{54}).
\qed
\ \\
By definition the spinor derivative on $S$ is given by the following
formula:\\
Let $\,\tilde{s}:U \longrightarrow Q_H\,$ be a local section in $Q_H$
and $\,s=(X_1,\ldots, X_{2n})=f_H(\tilde{s})\in P_H\,$ the
corresponding orthonormal basis in $(H,L_\theta)$. Consider a local spinor
field $\,\phi=[\,\tilde{s},u\,]\,$ in $S$. Then
\bean
\nabla^{S}\phi &=&
[\,\tilde{s},du\,]-\frac{1}{2}\,h(\nabla s_1,s_2)\,s_1\cdot s_2\cdot\phi
-\frac{1}{2}\suml^{2n}_{k=1}h(\nabla s_1,X^*_k)\,s_1\cdot
X^*_k\cdot\phi\\ && + \frac{1}{2}\suml^{2n}_{k=1}h(\nabla
s_2,X_k^*)\,s_2\cdot X^*_k\cdot\phi + \frac{1}{2}\suml_{k<l}
h(\nabla X_k^*,X^*_l)\,X^*_k\cdot X^*_l\cdot\phi.
\eean
Using the
definition of $s_1$ and $s_2$ (see (\ref{39})) and Proposition \ref{pr20} we obtain
$\,h(\nabla s_1, s_2)=0\,$. Furthermore, if we denote by $a_k(Z)$ the
vector field
\[ a_k(Z):= h(\nabla_Zs_2,X^*_k)\,s_2-h(\nabla_Zs_1,X^*_k)\,s_1, \]
from Proposition \ref{pr20} results
\bean
a_k(N) &=& 0\\
a_k(T^*) &=& -i\,\frac{c}{2}\,\Omega^A(T,X_k) \cdot N\\
a_k(X^*_j) &=& \frac{1}{2}\, d\theta
(X_j,X_k)\,T^*-\frac{1}{2} \{L_\theta([T,X_j], X_k) +\\
&& + L_\theta([T,X_k],X_j)+i\,\frac{c}{2}\,\Omega^A(X_j,X_k)\}N.
\eean
These formulas and Proposition \ref{pr20} give the following formulas
for the spinor derivative in the spinor bundle $S$ of $(B,h)$: \\

\begin{pr} \label{pr21}
Let $\,\tilde{s}:U \longrightarrow Q_H\,$ be a local section in
$Q_H$, $\,s=f_H(\tilde{s})=(X_1,\ldots,X_{2n})\,$ and let
$\,\phi=[\,\tilde{s},u\,]\,$ be a local section in $S$. Then for the
spinor derivative of $\phi$ holds:
\bean \nabla^{S}_N\phi &=&
[\,\tilde{s},N(u)\,]+\frac{1}{4} \,d\theta^*\cdot\phi\\
\nabla^{S}_{T^*}\phi &=& [\,\tilde{s},T^*(u)\,]+i\,\frac{c}{2}\,(T\ecke
\Omega^A)^*\cdot
N\cdot\phi-i\,\frac{c}{8}\,(\Omega^A_\theta)^*\cdot\phi\\
&& +
\frac{1}{4}\suml_{k<l}\{L_\theta([T,X_k],X_l)-L_\theta([T,X_l],X_k)
\}X^*_k\cdot X^*_l \cdot\phi\\
\nabla^{S}_{X^*}\phi &=&
[\,\tilde{s},X^*(u)\,] - \frac{1}{4}(X \ecke d\theta)^* \cdot T^*\cdot \phi
+ i\,\frac{c}{8}\,(X \ecke \Omega^A)_\theta^*\cdot N\cdot\phi\\ && +
\frac{1}{4}\suml^n_{k=1}\{L_\theta([T,X],X_k)+L_\theta([T,X_k],X)\}X^*_k\cdot
N\cdot\phi\\ && +
\frac{1}{2}\suml_{k<l}L_\theta(\nabla^W_X X_k,X_l)\,X^*_k\cdot X^*_l \cdot
\phi,
\eean
where $\sigma_\theta$ denotes the projection of a form
$\sigma\in \Lambda^p M$ onto $\Lambda^p_\theta M$, $\sigma^*_\theta$
is its horizontal lift on $B$ and the vector field $X$ belongs to the set
$\{X_1,\ldots,X_{2n}\}$.
\end{pr}
\ \\
\begin{pr} \label{pr22}
Let $(X_1,\ldots,X_{2n})$ be a local  orthonormal basis of
$(H,L_\theta)$ with $\,X_{2\alpha}=J(X_{2\alpha-1})\,$. Denote by
$\,\sigma^1, \ldots,\sigma^{2n}\,$ the dual basis of
$\,(X_1,\ldots,X_{2n})\,$
and by $\,s=(Z_1, \ldots,Z_n)\,$,
$\,Z_\alpha=\frac{1}{\sqrt{2}}(X_{2\alpha-1}-iX_{2\alpha})\,$, the
corresponding local unitary basis of $(T_{10},L_\theta)$. Consider the
2-forms
\begin{eqnarray*}
b_s &:= &
\suml_{k<l} \,\{L_\theta([T,X_k],X_l) - L_\theta([T,X_l],X_k)\}\,
\sigma^k \wedge\sigma^l, \\
 d_s(X) &:=& \suml_{k<l}\,L_\theta(\nabla^W_X
X_k,X_l)\,\sigma^k\wedge\sigma^l\,,\quad \qquad X\in H.
\end{eqnarray*}
Then
\begin{enumerate} \item[1)] $b_s\in\Lambda^{1,1}_\theta(M)\,$ and
$\,\mbox{Tr}_{\,\theta}\, b_s = 2 \mbox{Tr}\, \omega_s(T)$
\item[2)]
$d_s(X)\in\Lambda^{1,1}_\theta(M)\,$ and $\,\mbox{Tr}_{\,\theta}
\,d_s(X)= \mbox{Tr}\,\omega_s(X)$,
\end{enumerate}
where $\omega_s$ is the matrix
of connection forms of the Webster connection $\nabla^W$ with respect
to $s=(Z_1,\ldots,Z_n)$.
\end{pr}
{\bf Proof:} A 2-form $\sigma$ belongs to $\Lambda^{1,1}M$ iff
$\,\sigma(JX,JY)= \sigma(X,Y)\,$ for all $X,Y\in H$. From formula
(\ref{20}) of Proposition \ref{pr12} follows for
$X,Y\in\{X_1,\ldots,X_{2n}\}$
\bean b_s(JX,JY) &=& L_\theta([T,JX],JY)-L_\theta([T,JY],JX)\\
& \stackrel{(\ref{20})}{=} &
L_\theta([T,X],Y)-L_\theta ([T,Y],X)\\ &=& b_s(X,Y).
\eean
Therefore, $b_s\in\Lambda^{1,1}_\theta M$. Furthermore,
\bean \mbox{Tr}_\theta \,b_s
&=& i\suml^n_{\alpha=1}b_s(X_{2\alpha-1},X_{2\alpha})\\ &=&
i\suml^n_{\alpha=1}\{L_\theta([T,X_{2\alpha-1}],X_{2\alpha})-L_\theta
([T,X_{2\alpha}],X_{2\alpha-1})\}.
\eean
Inserting \[
X_{2\alpha-1}=\frac{1}{\sqrt{2}}(Z_\alpha+\bar{Z}_\alpha),\quad
X_{2\alpha} =\frac{i}{\sqrt{2}}(Z_\alpha-\bar{Z}_\alpha) \]
one obtains
\bean
\mbox{Tr}_\theta\, b_s &=&
\suml^n_{\alpha=1}\{L_\theta([T,Z_\alpha],Z_\alpha)-
L_\theta([T,\bar{Z}_\alpha],\bar{Z}_\alpha)\}\\ &=&
2i \,\suml^n_{\alpha=1}\,\mbox{Im}\,
\{L_\theta(\mbox{pr}_{10}[T,Z_\alpha], Z_\alpha)\}\\ &=&
2i\,\suml^n_{\alpha=1}\,\mbox{Im}\,L_\theta(\nabla^W_T
Z_\alpha,Z_\alpha))\\ &=&
2i\,\mbox{Im}\,(\,\mbox{Tr}\,\omega_s(T))\,.
\eean
Since $\nabla^W$ is metric
with respect to $L_\theta$, we have $\,\omega_{\alpha\beta}
+\overline{\omega_{\beta\alpha}}=0\,$. Hence,
$\,\omega_{\alpha\alpha}(T)\,$ is
imaginary. Therefore, $\,\mbox{Tr}_\theta\, b_s=2\,\mbox{Tr}\,
w_s(T)$.\\
According to formula (\ref{31}) of Proposition \ref{pr15} and formula
(\ref{19}) of Proposition \ref{pr12} we have for
$Y,Z\in\{X_1,\ldots,X_{2n}\}\,$ and
$\,X\in H$
\bean
d_s(X)(JY,JZ) &=& L_\theta(\nabla^W_XJY,JZ) \,=\,
L_\theta(J\nabla^W_XY,JZ)\\ &=& L_\theta(\nabla_X^WY,Z)\, = \,
d_s(X)(Y,Z).
\eean
This shows that $\,d_s(X)\in\Lambda^{1,1}_\theta(M)\,$. Furthermore,
\bean
\mbox{Tr}_\theta\, d_s(X) &=&
i\,\suml^n_{\alpha=1}\,L_\theta(\nabla^W_X X_{2\alpha-1}, X_{2\alpha})\\
&=&
\frac{1}{2}\,\suml^n_{\alpha=1}\,\{L_\theta(\nabla_X^WZ_\alpha,Z_\alpha)-
L_\theta(\nabla^W_X\bar{Z}_\alpha,\bar{Z}_\alpha)\}\\ &=& i\,\mbox{Im
Tr}\, \omega_s(X)\\ &=& \mbox{Tr}\, \omega_s(X).
\eean
\qed
\ \\
Next we proof a property of the spinor bundle $S_H$ of $(H,L_\theta)$,
which is very similar to the properties of the spinor bundle of
K\"ahler manifolds (see \cite{Kirchberg:86}).\\

\begin{pr} \label{pr23}
Let $\,(M^{2n+1},T_{10},\theta)\,$ be a strictly
pseudoconvex spin manifold and $(\sqrt{F},h_\theta)$ its Fefferman
space. Then the spinor bundle $S_H$ of $(H,L_\theta)$ has the
following properties:
\begin{enumerate}
\item
$S_H$ decomposes
into $\,n+1\,$ subbundles
\[ S_H=\bigoplus\limits^n_{r=0} S_{(-n+2r)i}, \]
where
$S_{ki}$ is the eigenspace of the endomorphism $d\theta\cdot :S_H \to
S_H$ to the eigenvalue $ki$. The dimension of $S_{ki}$ is
{\footnotesize $\left(
\begin{array}{c} n\\ \frac{n+k}{2}\end{array}\right)$}. In particular,
there are two 1-dimensional subbundles $\,S_{\varepsilon ni},\, \var =
\pm 1\,$, of $\,S_H\,$ satisfying $\,d\theta\cdot |_{S_{\var ni}}=\var
ni\cdot\mbox{Id}_{S_{\var ni}}$.
\item
If
$\,\sigma\in\Lambda^{1,1}_\theta M\,$, then \[ \sigma\cdot|_{S_{\var
ni}}=\var\cdot\mbox{Tr}_\theta(\sigma)\cdot\mbox{Id}_{S_{ \var ni}}.
\]
\item
The induced bundles $\pi^*S_{n\var i}$ on the Fefferman
space $\sqrt{F}$ are trivial. A global section $\psi_{\var}\in
\Gamma(\pi^*S_{n\var i})$ is given in the following way:\\
Let
$\,\tilde{s}:U \longrightarrow Q_H\,$ be a local section in $Q_H$, $s$
the local
unitary basis in $(T_{10},L_\theta)$, corresponding to
$\,f_H(\tilde{s}):U \longrightarrow P_H\,$. Furthermore, let
$\,\sqrt{\tau_s}:U \longrightarrow \sqrt{F}\,$ be the local section in
$\sqrt{F}$
defined by $\tilde{s}$ and let $\,\varphi_s:\sqrt{F}|_U
\longrightarrow S^1\,$ be given
by $\,p=\sqrt{\tau_s(\pi(p))}\cdot \varphi_s(p)\,$. Then \[
\psi_{\var}(p) :=
[\,\tilde{s}(\pi(p)),\varphi_s(p)^{-\var}u(\var,\cdots,\var)\,]
\] defines a global section in $\,\Gamma(\pi^*S_{n\var i})$.
\end{enumerate}
\end{pr}
{\bf Proof:} $\Delta^\pm_{2n,0}$ is a $U(n)$-representation, where
$U(n)$ acts by
\[ U(n)\,\stackrel{\ell}{\longrightarrow} \,
\mbox{Spin}^{\Bbb C}(2n)\,\, \;\stackrel{\Phi_{2n,0}}{\longrightarrow} \,\,\,
\mbox{GL}(\Delta^\pm_{2n,0}). \]
The element $\,\Omega_0=e_1\cdot
e_2+\cdots+ e_{2n-1}\cdot e_{2n}\in \mbox{Cliff}\,^{\Bbb C}_{2n,0}\,$
acts on
$\Delta^\pm_{2n,0}$ by \[ \Omega_0\cdot
u(\var_1,\ldots, \var_n) =i\,(\suml^n_{k=1}\var_k)\,u(\var_1,\ldots,
\var_n).
\]
Hence, $\Delta^\pm_{2n,0}$ decomposes into $U(n)$-invariant
eigenspaces $E_{\mu_r}(\Omega_0)$ of $\Omega_0$ to the eigenvalues
$\mu_r=(-n+2r)i$, $r=0,\ldots,n\,$. In particular, the eigenspace to the
eigenvalue $\var ni$, $\var=\pm 1\,,$ is 1-dimensional and given by \[
E_{in\var}(\Omega_0)=\C\cdot u(\var,\ldots,\var). \]
By definition of
$\ell$ (see (\ref{34})) we obtain for $\,A=\mbox{diag}(e^{i\theta_1},
\ldots e^{i\theta n)}$
\bea\label{55}
\ell(A)u(\var,\ldots,\var)=\left\{\begin{array}{cl}
u(\var,\ldots,\var) &\quad
\var=-1\\ \mbox{Det} A\cdot u(\var,\ldots,\var) &\quad \var=1.
\end{array}\right.
\eea
Hence, $E_{-ni}$ is the trivial
$U(n)$-representation and $E_{ni}$ is isomorphic to the
$U(n)$-representation $\Lambda^n(\C^n)$. Since the subspaces
$E_{\mu_r} (\Omega_0)$ of $\Delta^\pm_{2n,0}$ are invariant under the
action of $\lambda^{-1}(U(n))$ we obtain the decomposition \[
S_H= \bigoplus\limits^n_{r=0} S_{\mu_r}, \]
where $\,S_{\mu_r}:=
Q_H\times_{\lambda^{-1}(U(n))} E_{\mu_r}(\Omega_0)$.
\\ If
$\,\tilde{s}:U \longrightarrow Q_H\,$ is a local section in $Q_H$, $d\theta$
acts on $S_H$ by \[
d\theta\cdot [\,\tilde{s}\,,\,v\,]=[\,\tilde{s}\,,\,\Omega_0\cdot v\,]. \]
Therefore,
$S_{\mu_r}$ is the eigenspace of $d\theta\cdot$ to the eigenvalue
$\mu_r$.\\
Now, let $\,\eta=[\,q\,,\,u(\var,\ldots,\var)\,]\in S_{\var ni}\,$,
$\var=\pm 1$. Denote $\,f_H(q)=(X_1,\ldots,X_{2n})\in P_H\,$,
$\,X_{2\alpha}=JX_{2\alpha-1}\,$ and $\,s=(Z_1,\ldots,Z_n)\,$ the
corresponding unitary basis in $(T_{10},L_\theta)$ with
$\,Z_\alpha=\frac{1}{\sqrt{2}}(X_{2\alpha-1}-i JX_{2\alpha-1})\,$. Let
$\,(\theta^1,\ldots,\theta^n)\,$ be the dual basis of $(Z_1,\ldots,Z_n)$
and $(\sigma^1,\ldots,\sigma^{2n})$ the dual basis of
$(X_1,\ldots,X_{2n})$. If $\,\sigma\in\Lambda^{1,1}_\theta M\,$ is a
form
of type (1,1), then
\bean
\sigma &=&
\suml^n_{\alpha,\beta=1}\,\sigma_{\alpha\beta}\,\theta^\alpha\wedge
\overline{\theta^\beta}\\ &=& \frac{1}{2}
\suml_{\alpha\not=\beta}\sigma_{\alpha\beta}\,(\sigma^{2\alpha-1}
\wedge
\sigma^{2\beta-1}+\sigma^{2\alpha}\wedge\sigma^{2\beta})
+\frac{i}{2}\suml_{\alpha,\beta}\sigma_{\alpha\beta}\,(\sigma^{2\alpha}\wedge
\sigma^{2\beta-1}-\sigma^{2\alpha -1}\wedge\sigma^{2\beta}).
\eean
Hence,
\bean
\sigma\cdot\eta \,=\,
[\,q\,,&\frac{1}{2}\,\suml_{\alpha\not=\beta}\sigma_{\alpha\beta}
\,(e_{2\alpha-1}\cdot e_{2\beta-1}+e_{2\alpha}\cdot e_{2\beta})\cdot
u(\var,\ldots,\var)\\ & + \,
\frac{i}{2}\,\suml_{\alpha,\beta}\sigma_{\alpha\beta}\,(e_{2\alpha}
\cdot
e_{2\beta-1}-e_{2\alpha-1}\cdot e_{2\beta})\cdot u(\var,\ldots,\var)\,]
\eean
where $\,\sigma_{\alpha\beta}=\sigma(Z_\alpha,\bar{Z}_\beta)\,$.
Using formula (\ref{1}) we obtain
\bean
(e_{2\alpha-1}\cdot
e_{2\beta-1}+e_{2\alpha}\cdot e_{2\beta})\cdot u(\var, \ldots,\var)
&=& 0
\hspace{3cm} \alpha\not=\beta\\
(e_{2\alpha}\cdot
e_{2\beta-1}-e_{2\alpha-1}\cdot e_{2\beta}) \cdot u(\var,\ldots,
\var)&=& 0  \hspace{3cm}
\alpha\not=\beta\\
 (e_{2\alpha}\cdot
e_{2\alpha-1}-e_{2\alpha-1}\cdot e_{2\alpha}) \cdot u(\var,\ldots,
\var) &=&
-2\var i\, u(\var,\ldots,\var).
\eean
Therefore, \bean \sigma\cdot\eta &=&
[\,q\,,\,\var\,\suml_\alpha\, \sigma(Z_\alpha,\bar{Z}_\alpha)\cdot
u(\var,\ldots,\var)\,]\\
&=& \var\cdot \mbox{Tr}_\theta\,\sigma\cdot\eta.
\eean
Now, let us consider the section $\psi_{\var}\in\Gamma(\pi^* S_{\var
ni})$ defined by
\[ \psi_{\var}(p) := [\, \tilde{s}(\pi(p))\,,\,
\varphi_s(p)^{-\var}u(\var,\ldots, \var)\,]. \]
Let $\tilde{s},
\tilde{\hat{s}}:U \longrightarrow Q_H$ be two local sections,
$\tilde{s}=
\tilde{\hat{s}}\cdot g$ and let $h:U\longrightarrow S^1$ be the
function defined
by (\ref{35}):
\be\label{56} \ell(\lambda(g))\cdot
h=g,\quad h^2=\mbox{Det}(\lambda(g))^{-1}. \ee
Then
$\,\varphi_{\hat{s}}(p)=\varphi_s(p)\cdot h(\pi(p))\,$ and
\bean
\psi_{\var}(p) &=& [\,\tilde{\hat{s}}\cdot g\,,\,
\varphi_s(p)^{-\var}u(\var,\ldots, \var)\,]\\
&=&
[\,\tilde{\hat{s}}\,,\,\varphi_s(p)^{-\var}\, g\cdot
u(\var,\ldots,\var)\,]\\ &=&
[\,\tilde{\hat{s}}\,,\,\varphi_{\hat{s}}(p)^{-\var}\, h^{\var}\,
g\cdot u(\var,\ldots,\var)\,]\\
&\stackrel{(\ref{56})}{=}
& [\,\tilde{\hat{s}}\,,\, \varphi_{\hat{s}}
(p)^{-\var}h^{\var+1}\ell(\lambda(g))\,u(\var,\ldots,\var)\,]\\
&\stackrel{(\ref{55}),(\ref{56})}{=}&
[\,\tilde{\hat{s}}\,,\,\varphi_{
\hat{s}}(p)^{-\var}\, u(\var,\ldots,\var)\,].
\eean
Hence,
$\,\psi_{\var}\,$ is a global section in the bundle $\,\pi^*S_{\var
ni}\,$ on
$\sqrt{F}$. \qed
\ \\

\section{Twistor spinors on Fefferman spaces} \label{s7}
Let $\,(M^{2n+1},T_{10},\theta)\,$ be a strictly
pseudoconvex spin manifold, $(\sqrt{F},\pi,M)$ the square root of the
canonical $S^1$-bundle corresponding to the spinor structure of
$(M,g_\theta)$ and $h_\theta$ the Fefferman metric on $\sqrt{F}$.
Denote by $\psi_{\var}\in\Gamma(\pi^*S_H)$ the global sections in the
bundles $\pi^* S_{\var ni}$ over $\sqrt{F}$ defined in Proposition
\ref{pr23}.
Now, we are able to solve the twistor equation on the Lorentzian
spin manifold $\,(\sqrt{F},h_\theta)\,$.\\

\begin{th} \label{t1}
Let $\,S:=\pi^*S_H\oplus\pi^*S_H\,$ be the spinor bundle
of $\,(\sqrt{F},h_\theta)\,$. Then the spinor fields
$\,\phi_{\var}:=(\psi_{\var},0) \in\Gamma(S)\,$, $\var=\pm 1\,$, are
solutions of the twistor equation on $\,(\sqrt{F},h_\theta)\,$ with the
following properties:
\begin{enumerate}
\item
The canonical vector field $V_{\phi_{\var}}$ of $\phi_{\var}$ is
a regular isotropic Killing vector field.
\item
$V_{\phi_{\var}}\cdot \phi_{\var}=0\,.$
\item
$\nabla^S_{V_{\phi_{\var}}}\phi_{\var}=-\frac{1}{\sqrt{2}}\, \var\,
i\,\phi_{\var}\,.$
\item
$\|\phi_{\var}\|_\xi\equiv 1$.
\end{enumerate}
\end{th}
{\bf Remark:} If $n$ is even, then $\phi_1$ and $\phi_{-1}$ are
linearly
independent spinor fields in $S^+$. If $n$ is odd then $\phi_1\in
\Gamma(S^+)$ and $\phi_{-1}\in\Gamma(S^-)$ (see Proposition
\ref{pr19}). The
second property of Theorem \ref{t1} shows that $\phi_{\var}$ is a pure or partially pure
spinor field (see \cite{Trautman:94}).
A vector field is called
{\em regular}, if all of its integral curves are closed and of
the same shortest period.\\
\ \\
{\bf Proof of Theorem 1:} We use the formulas for the
spinor derivative in $S$ given in Proposition \ref{pr21} for the
Fefferman connection $\,A=A^{\sqrt{}}_\theta\,$ and the constant
$\,c=\frac{8}{n+2}\,$. Let $\tilde{s}:U\longrightarrow Q$ be a local
section and
$\,\varphi_s:\sqrt{F}|_U \longrightarrow S^1\,$ the corresponding
transition function
in $\sqrt{F}$ (see Proposition \ref{pr23}). Then for the fundamental
vector field $N$ on $\,\sqrt{F}\,$
\be\label{57}
N(\varphi_s)=\frac{n+2}{4}i \, \varphi_s
\ee
holds. If $\,Y^*\,$  is an
$\,A^{\sqrt{}}_\theta$-horizontal lift of a vector field $Y$ on $M$,
we obtain using standard formulas from connection theory
\bea\label{58}
Y^*(\varphi_s) &=& -\varphi_s\cdot
\sqrt{\tau_s}^*A^{\sqrt{}}_\theta(Y)\nonumber\\ &=&
\frac{1}{2}\,\varphi_s\,\{\,\mbox{Tr}\,
\omega_s(Y)+\frac{i}{2(n+1)}R^W
\theta(Y)\},
\eea
where $\omega_s$ is the matrix of connection forms of
the Webster connection with respect to the unitary basis $s$ in
$(T_{10},L_\theta)$ corresponding to $f_H(\tilde{s})$. According to
Proposition \ref{pr19} we have $\,N\cdot \phi_{\var}=0\,$.
Therefore, from
Proposition \ref{pr21} and (\ref{57}), (\ref{58}) result
\bean
\nabla^S_N\,\phi_{\var}&=&\left(-\var\,\frac{n+2}{4}
\,i\,\psi_{\var}+ \frac{1}{4}\,d\theta\cdot \psi_{\var}\,,\,0\,\right)\\
\nabla^S_{T^*}\phi_{\var}
&=&\left(-\frac{1}{2}\,\var\,\{\mbox{Tr}\,\omega_s(T)+\frac{i}{
2(n+1)}R^W\}\psi_{\var}+\frac{1}{4}\,b_s\cdot \psi_{\var} -
i\,\frac{1}{n+2}\,\Omega^{A^{\sqrt{}}_\theta}_\theta \cdot
\psi_{\var}\,,\,0\,\right)\\
\nabla^S_{X^*} \phi_{\var}&=& \left(-\frac{1}{2}\,\var\,\mbox{Tr}\,
\omega_s(X)\,\psi_{\var} +
\frac{1}{2}\,d_s(X)\cdot\psi_{\var}\,,\,0\,\right)
-\frac{1}{4}(X \ecke d\theta)^*\cdot T^*\cdot\phi_{\var}\,,
\eean
where
$\,b_s\,$ and $\,d_s(X)\,$ are the $\Lambda^{1,1}$-forms defined in
Proposition \ref{pr22}. Since $\psi_{\var}$ is a section in
$\pi^*S_{\var
ni}$, $b_s$ and $d_s(X)$ act on $\psi_{\var}$ by multiplication with
$\,\var\mbox{Tr}_\theta\, b_s\,$ and $\,\var\,\mbox{Tr}_\theta\,
d_s(X)\,$,
respectively (Proposition \ref{pr23}). Hence, according to Proposition
\ref{pr22},
\bean
\nabla^S_{T^*}\phi_{\var} &=&
\left(\,-i\frac{1}{n+2}\,\Omega^{A^{\sqrt{}}_\theta}_\theta
\cdot\psi_{\var} - \var\,\frac{i}{4(n+1)}\,R^W
\psi_{\var}\,,\,0\,\right)\\
\nabla^S_{X^*}\phi_{\var} &=& -\frac{1}{4}(X\ecke d\theta)^*\cdot
T^*\cdot \phi_{\var}.
\eean
Furthermore, $\psi_{\var}$ is an
eigenspinor of the action of $d\theta$ on $S_H$ to the eigenvalue
$\var ni$. Therefore,
\be\label{59}
\nabla^S_N\phi_{\var} =
-\frac{\var}{2}\,i\,\phi_{\var}\,.
\ee
Because of
\bean
\Omega^{A^{\sqrt{}}_\theta}_\theta \,=\,
-\frac{1}{2}\,\mbox{Ric}^W_\theta - \frac{i}{4(n+1)}
\,d(R^W\theta)_\theta\,=\,
-\frac{1}{2}\,\mbox{Ric}^W_\theta-\frac{i}{4(n+1)}\,R^Wd\theta,
\eean
the curvature $\,\Omega^{A^{\sqrt{}}_\theta}_\theta\,$ of the
Fefferman connection is a form of
type (1,1). Hence,
\bean
\Omega^{A^{\sqrt{}}_\theta}_\theta\cdot\psi_{\var} &=& \var \,\,
\mbox{Tr}_\theta (\Omega^{A^{\sqrt{}}_\theta}) \, \psi_{\var}\\ &=&
(-\frac{1}{2}\,\var \, R^W - \frac{i\var}{4(n+1)} \, R^W \,
in)\,\psi_{\var}\\ &=& -\var\frac{n+2}{4(n+1)} R^W \,\psi_{\var}.
\eean
Therefore, we obtain \be\label{60} \nabla^S_{T^*}\phi_{\var}=0. \ee
According to Proposition \ref{pr19},
$\,T^*\cdot\phi_{\var}=(\,0,\sqrt{2}\psi_{\var}\,)\,$. If
$X\in\{X_1,\ldots,X_{2n}\}$, the 1-form $X\ecke d\theta$ acts on the
spinor bundle by Clifford multiplication with $J(X)$. Hence, we have
\be\label{61}
\nabla^S_{X^*}\phi_{\var}\,=\,\left(\,0\,,\,-\frac{\sqrt{2}}{4}J(X)\cdot
\psi_{\var}\,\right).
\ee
Now, using $\,s_1= \frac{1}{\sqrt{2}}(N-T^*)\,$ , $\,
s_2=\frac{1}{\sqrt{2}}(N+T^*)\,$, we obtain
\bean
-s_1\cdot\nabla^S_{s_1}\phi_{\var} \,=\,
s_2\cdot\nabla_{s_2}^S \phi_{\var} \,=\, X^* \cdot\nabla_{X^*}^S
\phi_{\var} \,=\,
\left(\,0,-\frac{1}{2\sqrt{2}}\, \var\,i\,\psi_{\var}\,\right),
\eean
where
$X\in\{X_1,\ldots,X_{2n}\}\,$. This
shows, that $\,\phi_{\var}\,$ is a twistor spinor (see Proposition
\ref{pr1}).\\
From Proposition \ref{pr19} it follows
\bean
(\phi_{\var},\phi_{\var})_\xi
\,=\, \langle s_1 \cdot \phi_{\var}, \phi_{\var} \rangle \,=\, \langle
\,(0,-\psi_{\var}),(\psi_{\var},0)\,\rangle \,=\,
(\psi_\var,\psi_\var)_{S_H} \,=\,1\,.
\eean
Furthermore, we obtain for the canonical vector field
$\,V_{\phi_{\var}}\,$
\bean
V_{\phi_{\var}} &=& \langle
s_1\cdot\phi_{\var},\phi_{\var}\rangle\, s_1- \langle
s_2\cdot\phi_{\var},\phi_{\var}\rangle\, s_2 - \suml^{2n}_{k=1}
\,\langle X^*_k\cdot\phi_{\var},\phi_{\var}\rangle X^*_k\\ &=& s_1+s_2
\,=\, \sqrt{2}\, N.
\eean
Therefore, $V_{\phi_{\var}}$ is regular and
isotropic and satisfies $\,V_{\phi_{\var}} \cdot\phi_{\var}=0\,$.
Because of (\ref{59}) we have \[
\nabla^S_{V_{\phi_{\var}}}\phi_{\var}=-\frac{1}{\sqrt{2}}\,\var\,
i \, \phi_{\var}.
\]
It remains to show, that the vertical vector
field $N$ is a Killing vector field. This follows directly from the
formulas of Proposition \ref{pr20}:
\[
L_Nh_\theta(Y,Z)=h_\theta(\nabla_YN,Z)+h_\theta(Y,\nabla_ZN)=0
\]
for all vector fields $Y$ and $Z$ on $\,\sqrt{F}\,$. \qed
\ \\
Conversely, we have\\

\begin{th} \label{t2}
Let $(B^{2n+2},h)$ be a Lorentzian spin manifold and
let $\varphi\in\Gamma(S)$ be a nontrivial twistor spinor on $(B,h)$
such that
\begin{enumerate}
\item
The canonical vector field $V_\varphi$ of $\varphi$ is a regular
isotropic Killing vector field.
\item
$V_\varphi\cdot\varphi=0\,.$
\item
$\nabla^S_{V_\varphi}\varphi=i \,c\,\varphi\,, \qquad
c=\mbox{const} \in\R\backslash \{0\}$.
\end{enumerate}
Then $\,B\,$ is an
$\,S^1$-principal bundle over a strictly pseudoconvex spin
manifold\\
$(M^{2n+1},T_{10},\theta)\,$ and $\,(B,h)\,$ is locally isometric to the
Fefferman space $\,(\sqrt{F},h_\theta)$ of $(M,T_{10},\theta)$.
\end{th}
{\bf Proof:} Since $\,V_\varphi\,$ is regular, it defines an
$S^1$-action
on $B$
\bean
B\times S^1 & \longrightarrow & B\\ (p,e^{it}) &\longmapsto &
\gamma^V_{t\cdot\frac{L}{2\pi}} (p)
\eean
where $\,\gamma^V_t(p)\,$ is the
integral curve of $\,V=V_\varphi\,$ through $p$ and $L$ is the period
of the integral curves. Then $\,M:=B/_{S^1}\,$ is an
$\,2n+1\,$-dimensional
manifold and $V$ is the fundamental vector field defined by the
element $\frac{2\pi}{L}i$ of the Lie algebra $i\R$ of $S^1$ in the
$S^1$-principal bundle $(B,\pi,M;S^1)$. Now we use Sparling's
characterization of Fefferman spaces, proved by Graham in
\cite{Graham:87}.
Let $W$ denote the (4,0)-Weyl tensor, $C$ the (3,0)-Schouten-Weyl
tensor and $K$ the (2,0)-Schouten tensor of $(B,h)$. Graham proved:\\
If $V$ is an isotropic Killing vector field such that
\bea
V\ecke W &=& 0 \label{62} \\ V\ecke C &=& 0 \label{63}\\ K(V,V) &=&
\mbox{const}
<0, \label{64}
\eea
then there exists a pseudo-hermitian structure
$(T_{10},\theta)$ on $M$ such that $(B,h)$ is locally isometric to the
Fefferman space $(F,h_\theta)$ of $(M,T_{10},\theta)$. The local
isometry is given by $S^1$-equivariant bundle maps
$\,\phi_U:B|_U\longrightarrow
F|_U\,$.\\
We first prove that $\,V=V_\varphi\,$ satisfies
(\ref{62})-(\ref{64}). Property (\ref{63}) is valid for each twistor
spinor (see Proposition \ref{pr10}). Using $\,W(X\wedge
Y)\cdot\varphi=0\,$ (see
(\ref{11}) of Proposition \ref{pr5}) and the assumption
$\,V_\varphi\cdot\varphi=0\,$ we obtain
\bean
0 &=& \{W(X\wedge Y)\cdot
V-V\cdot W(X\wedge Y)\}\cdot\varphi\\
&=& 2\,\{ V\ecke \,W(X\wedge Y)\}\cdot
\varphi\\
&=& 2\, W(X,Y,V)\cdot\varphi
\eean
for all vector fields $X$ and $Y$ on $B$.
Since $V_\varphi$ is a nontrivial isotropic Killing field, it
has no zeros. Hence, by Proposition \ref{pr6}, the twistor spinor
$\varphi$
has no zeros and therefore, the vector field $W(X,Y,V)$ must be
isotropic for all vector fields $X,Y$ on $B$. Because of
\[
W(X,Y,V,V)=h(W(X,Y,V),V)=0\,, \]
$W(X,Y,V)$ is orthogonal to the
isotropic vector field $V$. Since $(B,h)$ has Lorentzian signature, it
follows that there is a 2-form $\lambda$ on $B$ such that
\be\label{65}
W(X,Y,V)=\lambda(X,Y)\, V \qquad\mbox{for all } X,Y\in \Gamma(TB).
\ee
Now, we use formula (\ref{12}) of Proposition \ref{pr5} to obtain
\bean
0 &=& V\cdot W(X\wedge Y)\cdot D\varphi-n \,\{V\cdot C(X,Y)+C(X,Y)\cdot
V \}\cdot \varphi\\
&=& V\cdot W(X\wedge Y)\cdot
D\varphi + 2n\,C(V,X,Y)\,\varphi.
\eean
Because of $V\ecke C=0$ it results
\be\label{66}
V\cdot W(X\wedge Y)\cdot D\varphi=0.
\ee
From the twistor equation (\ref{6}) and the assumption
$\,\nabla^S_V\varphi=
i\,c\,\varphi\,$ it follows
\bea\label{67}
W(X\wedge Y)\cdot V\cdot
D\varphi &=& -n\,\, W(X\wedge Y)\cdot \nabla^S_V\varphi\nonumber\\
&=& - n i c\,\,
W(X\wedge Y)\cdot\varphi\nonumber\\
&\stackrel{(\ref{11})}{=}& 0\,.
\eea
Then (\ref{65}), (\ref{66}) and
(\ref{67}) give
\bean
0 &=& W(X\wedge Y)\cdot V\cdot D\varphi-V\cdot
W(X\wedge Y)\cdot D\varphi\\
&=& 2\,\,W(X,Y,V)\cdot D\varphi\\ &=&
2\lambda(X,Y)\, V\cdot D\varphi\\
&\stackrel{(\ref{6})}{=}&
-2n\,\,\lambda(X,Y) \, \nabla_V^S
\varphi\\ &=& -2nci\,\lambda (X,Y)\,\varphi.
\eean
Therefore,
$\lambda\equiv 0$ and $V\ecke W=0$. Using formula (\ref{10}) of
Proposition \ref{pr5} we obtain
\bean
V\cdot \nabla_V^SD\varphi \,=\,
\frac{n}{2}\,\{V\cdot K(V)+K(V)\cdot V\}\cdot\varphi \,=\,
-n\,K(V,V)\varphi.
\eean
Since $V$ is an isotropic Killing field, it
satisfies $\nabla_VV=0$. It follows
\bean
\nabla_V^S(V\cdot D\varphi)
\,=\, \nabla_VV\cdot D\varphi+V\cdot\nabla_V^SD\varphi \,=\,
-n\,K(V,V)\,\varphi\,
\eean
and from the twistor equation
\[ \nabla^S_V\nabla^S_V\varphi = K(V,V)\,\varphi. \]
Using
$\,\nabla_V^S\varphi=ic\varphi\,$ we obtain $\,K(V,V)=-c^2\,$.
Therefore, the
canonical vector field $V_\varphi$ of the twistor spinor $\varphi$
satisfies the conditions of Sparling's characterization theorem for
Fefferman metrics. Now, we proceed as in Graham's proof of that
theorem. Since $\,V_{\alpha\varphi} =|\alpha|^2V_\varphi\,$ we can
normalize $\varphi$ in such a way that
$\,K(V_\varphi,V_\varphi)=-\frac{1}{4}\,$. Then, let $\tilde{T}$
be the vector field on $B$ defined by \[
h(\tilde{T},X)=-4\,K(X,V_\varphi) \,,\qquad X\in \Gamma(TB). \]
$\tilde{T}$
is isotropic and $\,h(\tilde{T},V_\varphi)=1\,$. Then we can use
$V_\varphi$ and $\tilde{T}$ to reduce the spinor structure of the
Lorentzian manifold $(B,h)$
to the group $\,\mbox{Spin}(2n)\,$. This reduced spinor structure projects
to a spinor structure of $(H,L_\theta)$, where $\theta$ is the
projection of the 1-form $\,\tilde{\theta}\in\Omega^1(B)\,$ dual to
$V_\varphi$ and $H\subset TM$ is the projection of the subbundle
$\,\tilde{H}=\mbox{span}(\tilde{T}, V_\varphi)^\bot\subset TB\,$ onto
$M$.
$J:H\to H$ is given by projection of the map
\bean
\tilde{J}: TB &\longrightarrow &
TB\\ X &\longmapsto & 2\,\nabla_XV_\varphi\,,
\eean
which acts on $\tilde{H}$
with $\tilde{J}^2=-id$. Then in \cite{Graham:87} is proved that
$\,(M,H,J,\theta)\,$ in fact is a strictly pseudoconvex manifold which
we equip with the spinor structure arising from that of $(H,L_\theta)$
by enlarging the structure group. In the same way as in
\cite{Graham:87} it
follows that $(B,h)$ is locally isometric to the Fefferman space
$\,(\sqrt{F},h_\theta)\,$, where the isometries are given by
$S^1$-bundle maps $\,\sqrt{F}|_U \longrightarrow B|_U$. \qed
\ \\
{\bf Remark:} Jerison and Lee studied the Yamabe problem on
CR-manifolds (see \cite{Jerison/Lee:87}). They proved that there is a
numerical CR-invariant $\lambda(M)$ associated with every
compact oriented strictly pseudoconvex manifold $M^{2n+1}$, which is
always less than or equal to the value corresponding to the sphere
$S^{2n+1}$ in $\C^n$ with its standard CR-structure. If $\lambda(M)$
is strictly less than $\lambda(S^{2n+1})$, then $M$ admits a
pseudo-hermitian structure $\theta$ with constant Webster scalar curvature
$\,R^W = \lambda(M)\,$.
Furthermore, one knows that the scalar curvature $R$ of the Fefferman
metric $h_\theta$ is a constant positive multiple of the lift of the
Webster scalar curvature $R^W$ to the Fefferman space (see
\cite{Lee:86}) . Now, let
$(M^{2n+1},T_{10})$ be a compact strictly pseudoconvex spin manifold
with $\,0 \not = \lambda(M) < \lambda(S^{2n+1})\,$. Choose a
pseudo-hermitian structure $\theta$ on
$(M,T_{10})$ such that the Webster scalar curvature $R^W$ is constant
(and non-zero since $\lambda(M) \not = 0$).\\
Let $\,\phi_{\var}$, $\var=\pm 1\,$, be the twistor
spinors on $(\sqrt{F},h_\theta)$, defined in Theorem \ref{t1}. Then
according to the remark following Proposition \ref{pr5} the spinor
fields
\[ \eta_{\var,\pm} \,:=\, \frac{1}{2}\, \phi_{\var}\, \pm \,
\sqrt{\frac{2n+1}{(2n+2)R}}\,\,\,D\phi_{\var} \]
are eigenspinors of the Dirac operator of the Lorentzian spin manifold
$\,(\sqrt{F},h_\theta)\,$ to the
eigenvalue $\,\pm \frac{1}{2} \sqrt{\frac{2n+2}{(2n+1)}R}\,$. The
length the spinor fields $\,\eta_{\var,\pm}\,$ is constant with respect
to the indefinite scalar product $\,\langle \cdot, \cdot \rangle\,$ as well as to
the positive definite scalar product $( \cdot, \cdot )_{\xi}$.

\bibliographystyle{alpha}
\bibliography{lit}
\end{document}